\documentclass[reqno,twoside,11pt,english]{amsart}
\usepackage{amsmath,amsfonts,amssymb,amsthm,epsfig}
\newcommand{\RNum}[1]{\uppercase\expandafter{\romannumeral #1\relax}}
\usepackage{comment}

\usepackage{color}
\usepackage[final,colorlinks,linkcolor=blue,anchorcolor=red,citecolor=blue]{hyperref}
\usepackage{graphicx}
\usepackage{titletoc,epsf}
\usepackage{amsmath,amsfonts,latexsym,amsthm,amsxtra,amssymb,bbm}
\allowdisplaybreaks
\usepackage{graphicx,float}
\usepackage{tikz}
\usepackage{ifthen}

\allowdisplaybreaks

\newtheorem{thm}{Theorem}[section]
\newtheorem{lem}[thm]{Lemma}
\newtheorem{op}[thm]{Open Problem}
\newtheorem{cor}[thm]{Corollary}
\newtheorem{prop}[thm]{Proposition}
\newtheorem{prob}[thm]{Problem}

\newtheorem{step}{Step}[section]

\newtheorem{cl}{Claim}[section]

\newtheorem{ca}{Case}
\newtheorem{sca}[section]{Subcase}
\newtheorem{scl}[section]{Subclaim}
\newtheorem{conj}[equation]{Conjecture}

\theoremstyle{definition}
\newtheorem{defn}[thm]{Definition}

\newtheorem{ques}[equation]{Question}
\newtheorem{rem}[thm]{Remark}
\newtheorem{exam}[thm]{Example}

\newcounter {own}
\def\theown {\thesection       .\arabic{own}}
\numberwithin{equation}{section}

\newenvironment{pf}[1][]{%
\vskip 3mm
\noindent
\ifthenelse{\equal{#1}{}}%
{{\slshape Proof. }}%
{{\slshape #1.} }%
}%
{\qed\bigskip}

\newtheorem{Thm}{Theorem}

\newtheorem{Ques}[Thm]{Question}
\newtheorem{Prop}[Thm]{Property}


\def\vint{\mathop{\mathchoice%
	{\setbox0\hbox{$\displaystyle\intop$}\kern 0.22\wd0%
		\vcenter{\hrule width 0.6\wd0}\kern -0.82\wd0}%
	{\setbox0\hbox{$\textstyle\intop$}\kern 0.2\wd0%
		\vcenter{\hrule width 0.6\wd0}\kern -0.8\wd0}%
	{\setbox0\hbox{$\scriptstyle\intop$}\kern 0.2\wd0%
		\vcenter{\hrule width 0.6\wd0}\kern -0.8\wd0}%
	{\setbox0\hbox{$\scriptscriptstyle\intop$}\kern 0.2\wd0%
		\vcenter{\hrule width 0.6\wd0}\kern -0.8\wd0}}%
\mathopen{}\int}

\newcommand{\IR}{{\mathbb R}}

\def\be{\begin{equation}}
\def\ee{\end{equation}}

\newcommand{\ben}{\begin{enumerate}}
\newcommand{\een}{\end{enumerate}}

\newcommand{\blem}{\begin{lem}}
\newcommand{\elem}{\end{lem}}
\newcommand{\bthm}{\begin{thm}}
\newcommand{\ethm}{\end{thm}}
\newcommand{\bcor}{\begin{cor}}
\newcommand{\ecor}{\end{cor}}
\newcommand{\beg}{\begin{exam}}
\newcommand{\eeg}{\end{exam}}

\newcommand{\bexs}{\begin{exs}}
\newcommand{\eexs}{\end{exs}}

\newcommand{\begs}{\begin{examples}}
\newcommand{\eegs}{\end{examples}}
\newcommand{\bdefe}{\begin{defn}}
\newcommand{\edefe}{\end{defn}}
\newcommand{\bprob}{\begin{prob}}
\newcommand{\eprob}{\end{prob}}
\newcommand{\bques}{\begin{ques}}
\newcommand{\eques}{\end{ques}}
\newcommand{\bei}{\begin{itemize}}
\newcommand{\eei}{\end{itemize}}
\newcommand{\bcon}{\begin{conj}}
\newcommand{\econ}{\end{conj}}
\newcommand{\bop}{\begin{op}}
\newcommand{\eop}{\end{op}}

\newcommand{\bas}{\begin{assertion}}
\newcommand{\eas}{\end{assertion}}

\newcommand{\bfa}{\begin{fact}}
\newcommand{\efa}{\end{fact}}

\newcommand{\bca}{\begin{ca}}
\newcommand{\eca}{\end{ca}}
\newcommand{\bsca}{\begin{sca}}
\newcommand{\esca}{\end{sca}}

\newcommand{\bcl}{\begin{cl}}
\newcommand{\ecl}{\end{cl}}

\newcommand{\bmlem}{\begin{mlem}}
\newcommand{\emlem}{\end{mlem}}

\newcommand{\bstep}{\begin{step}}
\newcommand{\estep}{\end{step}}

\newcommand{\bscl}{\begin{scl}}
\newcommand{\escl}{\end{scl}}

\newcommand{\bcons}{\begin{conjs}}
\newcommand{\econs}{\end{conjs}}

\newcommand{\bprop}{\begin{prop}}
\newcommand{\eprop}{\end{prop}}

\newcommand{\br}{\begin{rem}}
\newcommand{\er}{\end{rem}}
\newcommand{\brs}{\begin{rems}}
\newcommand{\ers}{\end{rems}}
\newcommand{\bo}{\begin{obser}}
\newcommand{\eo}{\end{obser}}
\newcommand{\bos}{\begin{obsers}}
\newcommand{\eos}{\end{obsers}}
\newcommand{\bpf}{\begin{pf}}
\newcommand{\epf}{\end{pf}}
\newcommand{\ba}{\begin{array}}
\newcommand{\ea}{\end{array}}
\newcommand{\beq}{\begin{eqnarray}}
\newcommand{\beqq}{\begin{eqnarray*}}
	\newcommand{\eeq}{\end{eqnarray}}
\newcommand{\eeqq}{\end{eqnarray*}}

\newcommand{\ds}{\displaystyle}

\newcounter{minutes}
\divide\time by 60
\newcounter{hours}
\multiply\time by 60 \addtocounter{minutes}{-\time}

\begin{document}

\bibliographystyle{amsplain}
\title[Dimension-free inner uniform estimates]{Dimension-free inner uniform estimates for quasigeodesics}
\author[C.-Y. Guo, M. Huang, Y. Li and X. Wang]{Chang-Yu Guo,  Manzi Huang$^*$, Yaxiang Li and Xiantao Wang}

\address[C.-Y. Guo]{Research Center for Mathematics and Interdisciplinary Sciences, Shandong University, 266237, Qingdao, P. R. China, and Department of Physics and Mathematics, University of Eastern Finland, 80101, Joensuu, Finland}
\email{{\tt changyu.guo@sdu.edu.cn}}

\address[M. Huang]{Manzi Huang, Department of Mathematics,
	Hunan Normal University, Changsha,  Hunan 410081, People's Republic
	of China} \email{mzhuang@hunnu.edu.cn}

\address[Y. Li]{Yaxiang Li, School of Mathematics and Statistics, Hunan First Normal University, Changsha,
	Hunan 410205, People's Republic of China} \email{yaxiangli@hnfnu.edu.cn; yaxiangli@163.com}

\address[X. Wang]{Xiantao Wang, Department of Mathematics,
	Hunan Normal University, Changsha,  Hunan 410081, People's Republic
	of China} \email{xtwang@hunnu.edu.cn}


\date{\today}
\keywords{Gromov hyperbolic, cone arc,
	quasigeodesic, inner uniform domain, John domian.\\
	${}^{\mathbf{*}}$ Corresponding author}

\begin{abstract}
	In this paper, establish a dimension-free inner uniform estimate for quasigeodesics.   
	More precisely, we prove that a $c_0$-quasigeodesic in a $\delta$-Gromov hyperbolic $c$-John domain in $\mathbb{R}^n$ is $b$-inner uniform, for some constant $b$ depending only on $c_0$, $\delta$ and $c$, but not on the dimension $n$. The proof relies crucially on the techniques introduced by Guo-Huang-Wang in their recent work [arXiv:2502.02930, 2025]. In particular, we actually show that the above result holds in general Banach spaces, which answers affirmatively an open question of J. V\"ais\"al\"a in [Analysis, 2004] and partially addresses the open question of Bonk-Heinonen-Koskela in [Asterisque, 2001]. As a byproduct of our main result, we obtain that  a $c_0$-quasigeodesic in  a $\delta$-Gromov hyperbolic $c$-John domain in $\mathbb{R}^n$ is a $b$-cone arc with a dimension-free constant $b=b(c_0,\delta,c)$. This resolves an open problem of J. Heinonen in [Rev. Math. Iberoam., 1989]. 
\end{abstract}

\date{\today}

\thanks{C.-Y. Guo is supported by the Young Scientist Program of the Ministry of Science and Technology of China (No.~2021YFA1002200), the NSF of China (No.~12101362 and 12311530037), the Taishan Scholar Project and the NSF of Shandong Province (No.~ZR2022YQ01). M. Huang and X. Wang are partly supported by NSF of China (No.12371071 and No.12571081). }


\maketitle 
\tableofcontents

\section{Introduction}

\subsection{Background}

The classical Riemann mapping theorem states that  a simply connected domain $D\subsetneq \IR^2$ is conformally equivalent to the unit disk $\mathbb{D}\subset \IR^2$. This celebrated theorem plays a fundamental role in planar conformal analysis and seeking for higher dimensional extension or even more abstract metric space extension has been one of the central problems in modern analysis and geometry. In a seminal work \cite{BHK}, Bonk, Heinonen and Koskela have successfully established a suitable extension of this uniformization theorem for the class of Gromov hyperbolic spaces. Beside a rather rich uniformization theory, they further 
\begin{itemize}
	\item studied the relationship between Gromov hyperbolicity and inner uniformality and established a dimension-dependent characterization between Gromov hyperbolicity and inner uniformality; see \cite[Theorem 1.11 and Proposition 7.12]{BHK}.
	\item proved that Gromov hyperbolic domains in $\IR^n$ satisfy the \emph{the Gehring-Hayman inequality} and a \emph{ball-separation condition}, quantitatively; see \cite[Section 7]{BHK}.
\end{itemize}
Moreover, motivated by the theory of quasiconformal mappings in infinite dimensional Banach spaces, they asked the following interesting challenging open problem  \cite[Page 5]{BHK}:
\medskip 

\textbf{Question (Bonk-Heinonen-Koskela):} \emph{Is there a general relationship between uniformity and Gromov hyperbolicity that would cover infinite-dimensional situations as well.}
\medskip 

As was pointed out in \cite[Page 1]{BHK}, ideas around Gehring-Hayman inequality play a crucial role in the work of Bonk-Heinonen-Koskela, which relies heavily on finite dimensional techniques such as Lebesgue measure and integration. To solve the above open question of Bonk-Heinonen-Koskela, we developed some new techniques to obtain dimension-free Gehring-Hayman inequality in our previous work \cite{Guo-H-W-2025}. In this paper, we continue our investigation on the relationship between inner uniformality and Gromov hyperbolicity in infinite dimensional Banach spaces. Adapting the techniques introduced in \cite{Guo-H-W-2025} and developing some new estimates in the setting of Gromov hyperbolic John domains, we shall improve the corresponding results of Bonk-Heinonen-Koskela \cite[Theorem 1.11 and Proposition 7.12]{BHK} into dimension-free form and thus it holds in infinite dimensional Banach spaces. 

Our main focus of this paper is to study the geometry of quasigeodesics in Gromov hyperbolic John domains. For this, we first recall a couple of relevant definitions. 
\begin{defn}\label{def:uniform domain}
	A domain $D$ in a metric space $X=(X,d)$ is called  {\it $c$-uniform}, $c\geq 1$, if each pair of points $z_{1},z_{2}$ in $D$ can be joined by a rectifiable curve $\gamma$ in $D$ satisfying
	\begin{enumerate}
		\item\label{con1} $\ds\min_{j=1,2}\{\ell (\gamma [z_j, z])\}\leq c\, d_D(z)$ for all $z\in \gamma$, and
		\item\label{con2} $\ell(\gamma)\leq c\,d(z_{1},z_{2})$,
	\end{enumerate}
	where $\ell(\gamma)$ denotes the arc-length of $\gamma$ with respect to the metric $d$,
	$\gamma[z_{j},z]$ the part of $\gamma$ between $z_{j}$ and $z$, and $d_D(z):=d(z,\partial D)$. In a $c$-uniform domain $D$, any curve $\gamma\subset D$, which satisfies conditions (1) and (2) above, is called a {\it $c$-uniform curve} or a {\it double $c$-cone curve}.
\end{defn}

If the condition \eqref{con2} in Definition \ref{def:uniform domain} is replaced by the weaker inequality 
\begin{equation}\label{eq:def for inner quasiconvexity}
	\ell(\gamma)\leq c\,\sigma_D( z_{1}, z_{2}),
\end{equation}
where $\sigma_D$ is the {\it inner distance} defined by
$$\sigma_D(z_1,z_2)=\inf \{\ell(\alpha):\; \alpha\subset D\;
\mbox{is a rectifiable curve joining}\; z_1\; \mbox{and}\; z_2 \},$$
then $D$ is said to be $c$-{\it inner uniform} and the corresponding curve $\gamma$ is called a $c$-\emph{inner uniform curve}.

If $D$ only satisfies the condition \eqref{con1} in Definition \ref{def:uniform domain},  then it is said to be a {\it $c$-John domain}, and the corresponding curve $\gamma$ is called a \emph{$c$-John curve} or a {\it $c$-cone curve}.

It follows immediately from the definition that 
\[
\{c\text{-uniform domains}\} \subset \{c\text{-inner uniform domains}\}\subset  \{c\text{-John domains}\}.
\]
It is not difficult to construct examples to see that the above inclusions are strict. For instance, a slit disk, $\mathbb{D}$ with irrational points on $[0,1)$ removed, in $\mathbb{R}^2$ is John but not inner uniform.  

The class of John domains was initially introduced by F. John in his study of elasticity \cite{Jo} and the name was coined by Martio and Sarvas in \cite{MS}, where they also introduced the class of uniform domains. These classes of domains are central in modern geometric function theory in $\mathbb{R}^n$ or more general metric spaces and have wide connections with many other mathematical subjects related to analysis and geometry; see for instance \cite{GM,Martin-1985,Jo81,H,HeKo,Hajlasz-Koskela-2000,Hei-Book-2001,Guo-Koskela-2017,HLPW}.

Another useful class of domains, which is larger than the class of (inner) uniform domains, is the so-called Gromov hyperbolic domains. Recall that a geodesic metric space $X=(X,d)$ is called $\delta$-Gromov hyperbolic, $\delta>0$, if each side of a geodesic triangle in $X$ lies in the $\delta$-neighborhood of the other two sides. This notion was introduced by M. Gromov in his celebrated work \cite{Gr-1} and it generalizes the metric properties of classical hyperbolic geometry.  Gromov hyperbolicity is a large-scale property and has been widely studied in geometric function theory and metric geometry; see for instance \cite{BB2,BB4,BS,BHK,Klm2014,Guo15,ZLR,HRWZ}.
\begin{defn}\label{def:Gromov hyperbolic domain}
	An Euclidean domain $D\subset \mathbb{R}^n$ is called $\delta$-Gromov hyperbolic if the metric space $(D,k_D)$ is $\delta$-Gromov hyperbolic, where $k_D$ is the quasihyperbolic metric in $D$.  
\end{defn}

There are no inclusion relation between Gromov hyperbolic domains and John domains: 
\begin{itemize}
	\item Gromov hyperbolic domains are not necessarily inner uniform or John. For instance, a simply connected planar domain with an outward pointing cusp is Gromov hyperbolic, but not inner uniform nor John.
	\item  John domains are not necessarily Gromov hyperbolic, for example, $\mathbb{R}^2\setminus \{(n,0): n\in \mathbb{N}\}$ is a John domain but not Gromov hyperbolic. 
\end{itemize}

A central tool in the study of various relationships between these different classes of domains is the quasihyperbolic geodesic (see Section \ref{subsec:quasihyperbolic geodesic} below for precise meaning), which is a natural extension of planar hyperbolic geodesic in higher dimensions. For instance, Gehring and Osgood \cite{GO} proved that each  quasihyperbolic geodesic in a $c$-uniform domain in $\mathbb{R}^n$ is a $b$-uniform arc with $b$ depending only on $c$.


The following result, due to Bonk, Heinonen and Koskela \cite{BHK}, gives a useful characterization of inner uniform domains and Gromov hyperbolic John domains in $\mathbb{R}^n$. 
\begin{Thm}[{\cite[Theorem 1.11 and Proposition 7.12]{BHK}}]\label{thm:A}
	Let $D\subsetneq \mathbb{R}^n$ be a domain with $n\geq 2$. Then the following conclusions hold. 
	\begin{enumerate}
		\item
		If $D$ is a $b$-inner uniform domain, then it is $\delta$-Gromov hyperbolic $b$-John with $\delta=\delta(b)$.
		\item
		If $D$ is a $\delta$-Gromov hyperbolic $c$-John domain, then each quasihyperbolic geodesic in $D$ is $b$-inner uniform with $b=b(n,\delta,c)$. In particular, $D$ is a $b$-inner uniform domain. 
	\end{enumerate}
\end{Thm}



Note that the first statement in Theorem \ref{thm:A} is dimension-free, while the second statement is dimension-dependent. It is thus natural to ask the following question. 

\begin{Ques}\label{Quest:main question}
	If $D\subset \mathbb{R}^n$ is a $\delta$-Gromov hyperbolic $c$-John domain, then is it true that each quasihyperbolic geodesic in $D$ is $b$-inner uniform with a dimension-free constant $b=b(\delta,c)$?  
\end{Ques}

Observe that an affirmative answer to Question \ref{Quest:main question} would most likely provide a solution to the following open question of  J. V\"ais\"al\"a \cite[Open question 3.17]{Vai2004}: \emph{Is every Gromov hyperbolic John domain in $E$ (or more generally, in an arbitrary Banach space) inner uniform}? 


Our main result of this paper will provide affirmative solutions to Question \ref{Quest:main question} as well as the above open problem of V\"ais\"al\"a. We will actually consider Question \ref{Quest:main question} (and V\"ais\"al\"a's open question) for the more general class of quasigeodesics. 


%



\subsection{Main results}
To state our main result, we first recall the following definition of quasigeodesics introduced by J. V\"ais\"al\"a \cite{Vai4}. 
\begin{defn}\label{def:quasigeodesic}
	A curve $\gamma\subset D$ is said to be a {\it
		$c$-quasigeodesic}, $c\geq 1$, if for all $x$, $y$ in $\gamma$, it holds
	$$\ell_{k}(\gamma[x,y])\leq c\, k_D(x,y),
	$$ where $\ell_{k}(\gamma[x,y])$ denotes the quasihyperbolic length of $\gamma[x,y]$ (see Section \ref{sec-2} for the definition).
\end{defn}

Our main result of this paper gives a dimension-free inner uniform estimate in Theorem \ref{thm:A}(2). 
\begin{thm}\label{thm:main Euclidean}
	Let $D\subsetneq \mathbb{R}^n$ be a $\delta$-Gromov hyperbolic $c$-John domain. Then every $c_0$-quasigeodesic $\gamma\subset D$ is $b$-inner uniform with $b=b(\delta,c,c_0)$.  
\end{thm}

In the special case $c_0=1$, Theorem \ref{thm:A} provides an affirmative solution to Question \ref{Quest:main question}. As a direct application of this result, we get an improved dimension-free version of Theorem \ref{thm:A}. 
\begin{thm}\label{thm:main 2 Euclidean}Let $D\subsetneq \mathbb{R}^n$ be a domain. Then the following conclusions hold.
	\begin{enumerate}
		\item
		If $D$ is a $b$-inner uniform domain, then it is $\delta$-Gromov hyperbolic $b$-John with $\delta=\delta(b)$.
		\item
		If $D$ is a $\delta$-Gromov hyperbolic $c$-John domain, then each $c_0$-quasigeodesic in $D$ is $b$-inner uniform with $b=b(\delta,c,c_0)$. In particular, $D$ is a $b$-inner uniform domain. 
	\end{enumerate}
\end{thm}

For the proof of Theorem \ref{thm:main Euclidean}, we shall employ the new techniques introduced in our previous work \cite{Guo-H-W-2025}, together with some nontrivial adaptions. As an application of this general approach, we shall indeed prove a version of Theorem \ref{thm:main Euclidean} in the setting of general Banach spaces, from which Theorem \ref{thm:main Euclidean} follows as a special case. 

\begin{thm}\label{thm:main}
	Let $E$ be a Banach space and $D\subsetneq E$ a $\delta$-Gromov hyperbolic $c$-John domain. Then every $c_0$-quasigeodesic $\gamma\subset D$ is $b$-inner uniform with $b=b(\delta,c,c_0)$. In particular, $D$ is $b$-inner uniform.  
\end{thm}

Notice that Theorem \ref{thm:main} provides an affirmative solution to the aforementioned open problem of J. V\"ais\"al\"a \cite[Open question 3.17]{Vai2004}. It also gives a partial affirmative answer to  \textbf{Question of Bonk-Heinonen-Koskela} by adding an extra John assumption along with Gromov hyperbolicity. On the other hand, by \cite[Remark 3.16]{Vai2004}, a bounded Gromov hyperbolic LLC-2 domain in a Banach space $E$ might not be inner uniform. This shows that the John domain assumption in Theorem \ref{thm:main} is quite sharp. 


To state our next application, we need to recall the definition of coarsely quasihyperbolic homeomorphism.
\begin{defn}\label{def:CQH map}
	A homeomorphism $f:D\to D'$, where $D'\subset E'$, is called
	{\it $C$-coarsely $M$-quasihyperbolic}, abbreviated $(M,C)$-CQH,
	if for all $x$, $y\in D$, it holds
	$$\frac{1}{M}(k_D(x,y)-C)\leq k_{D'}(f(x),f(y))\leq M\;k_D(x,y)+C,$$
	where $k_D$ and  $k_{D'}$ denote the quasihyperbolic metrics in $D$ and $D'$, respectively (see Section \ref{sec-2} for the definition). In particular, an $(M,0)$-CQH mapping is said to be
	{\it $M$-bi-Lipschitz} in the quasihyperbolic metric.
\end{defn}

By \cite[Theorem 3]{GO}, every $K$-quasiconformal mapping $f:\Omega\to \Omega'\subset \IR^n$ is an $(M, C)$-CQH homeomorphism, where the constants $M$ and $C$ depend only on $K$ and $n$. The reverse implication does not hold in general; see \cite[Example 2.1]{Guo-H-W-2025}. As Gromov hyperbolicity is preserved by CQH maps \cite[Theorem 3.18]{Vai10}, we obtain from Theorem \ref{thm:main} the following corollary.

\begin{cor}\label{coro1}
	Suppose that $E$ is a Banach space and $D\subset E$ is a $c$-John domain, which
	homeomorphic to a $\delta$-Gromov hyperbolic domain via a $(M,C)$-CQH mapping. Then each $c_0$-quasigeodesic in $D$ is $b$-inner uniform with $b=b(\delta,c,M,C,C_0)$.
\end{cor}

\subsection{Short comments on the proof}
We now briefly comment on the proof of Theorem \ref{thm:main}. It should be noticed that the original approach of Bonk, Heinonen and Koskela \cite{BHK} relies crucially on the conformal $n$-modulus of path families and Ahlfors $n$-regularity of $n$-dimensional Lebesgue meassure of $\mathbb{R}^n$, and thus it cannot be used in the setting of infinite dimensional Banach spaces. 

To obtain a dimension-free estimate on the inner uniformality coefficient of quasigeodesics, one needs to verify the dimension-free cone property in Definition \ref{def:uniform domain}(1) and the dimension-free Gehring-Hayman inequality \eqref{eq:def for inner quasiconvexity} for quasigeodesics. 

In \cite{ZLR}, the authors obtained a dimension-free cone property for quasihyperbolic geodesics in locally compact geodesic spaces. However, their proof uses the uniformization theory of Bonk, Heinonen, and Koskela \cite{BHK} and thus does not extend to the setting of Banach spaces and does not work for the more general class of quasigeodesics. In this paper, we shall provide a more direct approach to obtain dimension-free cone property for quasigeodesics (see Theorem \ref{thm-1}). In particular, this resolves an open question of J. Heinonen \cite[Question 2]{H}.

The more difficult part of the proof of Theorem \ref{thm:main} lies in the second step, that is, to establish the dimension-free Gehring-Hayman inequality for quasigeodesics. This relies crucially on the new techniques introduced in the recent work of Guo, Huang and Wang \cite{Guo-H-W-2025}, where the authors have successfully obtained a dimension-free Gehring-Hayman inequality for quasigeodesics in domains that are coarsely quasihyperbolically equivalent to uniform domains in Banach spaces. The difference between these two papers is that we have relaxed the CQH uniformality of \cite{Guo-H-W-2025} to the weaker Gromov hyperbolicity in this paper. On the other hand, as a compensation, we have the extra John domain assumption at hand. \emph{We do not know whether the John domain assuption is really necessary for this step or not}. Comparing with the proof in \cite{Guo-H-W-2025}, there are a few nontrivial adaptions in this paper: 
\begin{itemize}
	\item develop some new estimates on the quasihyperbolic distances of points on quasigeodesics (such as Lemmas \ref{thm3} and \ref{lem02});
	
	\item establish new estimates in the construction of a sequence of new six-tuples (see Lemmas \ref{lem-22-09-30}  \ref{lem-24-05});
	
	\item find a new contradiction argument (see Section \ref{subsec:final section}).
\end{itemize}

{\bf Structure.}
The organization of this paper is as follows. Section \ref{sec-0} contains preliminaries and auxiliary results about quasihyperbolic metric, quasigeodesics and Gromov hyperbolic spaces. In section \ref{sec-2}, we will prove the dimension-free cone property for quasigeodesics in Theorem \ref{thm-1}. In section \ref{sec-3}, we prove the dimension-free Gehring-Hayman inequality and thus completes the proof of Theorem \ref{thm:main}.   

{\bf Notation.} Throughout the paper, we always assume that $E$ denotes a real
Banach space with dimension at least $2$. The norm of a vector $z$
in $E$ is written as $|z|$, and for each pair of points $z_1$, $z_2$
in $E$, the distance between them is denoted by $|z_1-z_2|$, the
closed line segment with endpoints $z_1$ and $z_2$ by $[z_1, z_2]$.
We always use $\mathbb{B}(x_0,r)$ to denote the open ball $\{x\in
E:\,|x-x_0|<r\}$ centered at $x_0$ with radius $r>0$. Similarly, for
the closed balls and spheres, we use the usual notations
$\overline{\mathbb{B}}(x_0,r)$ and $ S(x_0,r)$,
respectively. 

\section{Preliminaries and auxiliary results}\label{sec-0}

\subsection{List of coefficients} \label{subsec:list of coefficients}
\begin{enumerate}
	\item  $\min\{c,C_0,M_0,\mu_0\}\geq 16$, $A_0=3C_0+2M_0$;
	
	\item  $c_1\geq 8c_0cA_0M_0$, where $c_1$ is from Lemma \ref{thm3};
	
	\item $c_2=e^{20c_0^2c_1^2A_1}$, where $A_1=A_1(\delta,c_0)>2c_0c_1e^{\delta}$ (cf. Lemma \ref{lem24-1}) and Theorem \ref{thm-1};
	
	
	\item $C_1=\ell_k(\beta_0)$, $\kappa=\kappa_1^4$ and $\kappa_1=e^{(96c_1c_2A_0)^6}$;
	
	\item $\tau_1=11A_0$, $\tau_2=\frac{1}{64}\log \kappa$ and $\tau_3=\frac{1}{16}\log \kappa$.
\end{enumerate}

\subsection{Quasihyperbolic distance and quasigeodesic}\label{subsec:quasihyperbolic geodesic} 
We start this section by recalling the definition of quasihyperbolic metric, which was initially introduced by Gehring and Palka \cite{GP} for a domain in $\IR^n$ and then has been extensively studied in \cite{GO}. Since we are concerning it in a Banach space, we shall follow the presentation from \cite{Vai3}. 

The {\it quasihyperbolic length} of a rectifiable arc 
$\gamma$ 
in a proper domain $D\subsetneq E$ is defined as
$$
\ell_{k}(\gamma):=\int_{\gamma}\frac{|dz|}{d_D(z)}.
$$
For any $z_1$, $z_2$ in $D$, the {\it quasihyperbolic distance}
$k_D(z_1,z_2)$ between $z_1$ and $z_2$ is set to be
$$k_D(z_1,z_2)=\inf_{\gamma}\{\ell_k(\gamma)\},
$$
where the infimum is taken over all rectifiable arcs $\gamma$
joining $z_1$ and $z_2$ in $D$.

An arc $\gamma$ from $z_1$ to $z_2$ is called a {\it quasihyperbolic geodesic} if
$\ell_k(\gamma)=k_D(z_1,z_2)$. Clearly, each subarc of a quasihyperbolic
geodesic is a quasihyperbolic geodesic. It is a well-known fact that a
quasihyperbolic geodesic between any two points in $E$ exists if the
dimension of $E$ is finite; see \cite[Lemma 1]{GO}. This is not
true in infinite dimensional spaces; see \cite[Example 2.9]{Vai4}. In order to remedy this shortage, J. V\"ais\"al\"a \cite{Vai4}
introduced the class of \textit{quasigeodesics} in $E$; see Definition \ref{def:quasigeodesic}. Moreover, he proved in \cite[Theorem 3.3]{Vai4} that
for each pair of points $z_1$ and $z_2$ in $D$ and each $c>1$, there exists a
$c$-quasigeodesic in $D$ joining $z_1$ and $z_2$.  Meanwhile, it was proved in \cite[Theorems 9.6 and 10.9]{Vai8} that a $c$-quasigeodesic in an $A$-uniform domain is an $A_1$-uniform arc with $A_1 = A_1(A, c)$.

For any $z_1$, $z_2$ in $D$, let $\gamma\subset D$ be an arc with end points $z_1$ and $z_2$. Then we have the following elementary estimates (see for instance \cite[Section 2]{Vai3}):
\beq\label{(2.1)}
\ell_{k}(\gamma)\geq
\log\Big(1+\frac{\ell(\gamma)}{\min\{d_D(z_1), d_D(z_2)\}}\Big)
\eeq
and
\beq\label{(2.2)}
\begin{aligned}
	k_{D}(z_1, z_2) &\geq  \log\Big(1+\frac{\sigma_D(z_1,z_2)}{\min\{d_D(z_1), d_D(z_2)\}}\Big)
	\\ 
	&\geq 
	\log\Big(1+\frac{|z_1-z_2|}{\min\{d_D(z_1), d_D(z_2)\}}\Big)
	\geq
	\Big|\log \frac{d_D(z_2)}{d_D(z_1)}\Big|.
\end{aligned}
\eeq


We have the following elementary lemma. 
\begin{lem}\label{lem-3-1}
	Suppose $D\subset E$ is a domain, $u$, $v\in D$, and $\alpha$ is a rectifiable curve in $D$ joining $u$ and $v$. Then, for each $w\in\alpha$ and $c\geq 1$, the following conclusions hold.
	
	$(1)$ If $\ell(\alpha[u,w])\leq c d_D(w)$, then
	$$\ds d_D(w)\geq \max \left\{\frac{2\ell(\alpha[u,w])+d_D(u)}{4c},\frac{1}{2c}d_D(u)\right\}.
	$$
	
	$(2)$  If $\ell(\alpha[u,w])\leq c d_D(w),$
	then
	$$k(u,w)\leq 4c \log\Big(1+\frac{\ell(\alpha[u,w])}{d_D(u)}\Big).$$

	$(3)$ If $|u-w|\leq c d_D(w),$  then
	$$\ds d_D(w)\geq \frac{1}{2c}d_D(u).
	$$
\end{lem}
\bpf  
(1) 
If $\ell(\alpha[u,w])\geq \frac{1}{2}d_{D}(u)$, then $$d_D(u)\leq 2\alpha(\beta[u,w])\leq 2c d_D(w)$$
and so 
\[
\max \left\{\frac{2\ell(\alpha[u,w])+d_D(u)}{4c},\frac{1}{2c}d_D(u)\right\}=\frac{2\ell(\alpha[u,w])+d_D(u)}{4c}\leq d_D(w). 
\] 

If $\ell(\alpha[u,w])< \frac{1}{2}d_{D}(u)$, then $$d_D(w)\geq d_D(u)-\ell(\alpha[u,w])>\frac{1}{2}d_D(u).$$
Since $c\geq 1$, we infer from the above estimate that
$$d_D(w)\geq \frac{1}{2c}d_D(u)=\max \left\{\frac{2\ell(\alpha[u,w])+d_D(u)}{4c},\frac{1}{2c}d_D(u)\right\}.$$
This proves (1). 

(2) It follows from (1) of this lemma that \beq k_{\Omega}(u,w)\nonumber&\leq&\ell_k(\alpha[u,w])=\int_{x\in\alpha[u,w]}\frac{|dx|}{d_{\Omega}(x)}\leq 4c\int_{x\in\alpha[u,w]}\frac{|dx|}{2\ell(\alpha[u,x])+d_D(u)}\\ \nonumber&\leq& 4c\log\Big(1+\frac{\ell(\alpha[u,w]}{d_{\Omega}(u)}\Big).\eeq

(3) This can be proved similarly as (1) and thus is omitted here. 
\epf

Next we discuss the coarse length of a curve intrudced by J. V\"ais\"al\"a in \cite[Section 4]{Vai4}. Let $\gamma$ be an (open, closed or half open) arc in $D$. Write $\overline{x}=(x_0,\ldots,x_m)$ with
$m\geq 1$, whose coordinates form a finite sequence of $m$ successive points of $\gamma$.
For $h\geq 0$, we say that $\overline{x}$ is {\it $h$-coarse} if
$k_D(x_{j-1},x_j)\geq h$ for all $j\in\{ 1, \ldots, m\}$. Let $\Phi_k(\gamma,h)$
be the family of all $h$-coarse sequences of $\gamma$. Set
$$s_k(\overline{x}):=\sum^{m}_{j=1}k_D(x_{j-1},x_j)\;\; \mbox{and}\;\;
\ell_k(\gamma, h):=\sup \left\{s_k(\overline{x}):\; \overline{x}\in \Phi_k(\gamma,h)\right\}$$
with the convention that $\ell_k(\gamma, h)=0$ if
$\Phi_k(\gamma,h)=\emptyset$. Then the number $\ell_k(\gamma, h)$ is called the
{\it $h$-coarse quasihyperbolic length} of $\gamma$.

\bdefe[Solid arc] \label{def:h-solid} An arc
$\gamma$ in $D$ is called {\it $(c,h)$-solid} with $c\geq 1$ and
$h\geq 0$ if
$$\ell_k(\gamma[z_1,z_2], h)\leq c\;k_D(z_1,z_2)$$ for all pairs of points $z_1$ and $z_2$ in $\gamma$. 
{We denote the class of all $(c,h)$-solid arcs in $D$ by $S^{c,h}(D)$. The symbol $S^{c,h}_{xy}(D)$ stands for those arcs in $S^{c,h}(D)$ with end points $x,y$. }
\edefe

By \cite[Section 2.2]{Vai7} or \cite[Section 4]{Vai4}, we know that for any arc $\gamma \subset D$, $\ell_k(\gamma, 0)=\ell_k(\gamma)$.
Notice that $\gamma$ is a $(c,0)$-solid arc if and only if it is a $c$-quasigeodesic.

\bdefe[Short arc] \label{def:h-short}
A curve $\gamma$ in $D$ is said to be {\it $h$-short}
with $h\geq 0$ if for all $x,y\in\gamma$, $$\ell_{k}(\gamma[x,y])\leq  k_D(x,y)+h.$$
\edefe

It is clear from Definition \ref{def:h-short} that every subarc of an $h$-short arc is $h$-short, every $h$-short arc must be $h_1$-short for any $h_1\geq h$, and every quasihyperbolic geodesic is $h$-short with $h=0$.

\bdefe[$\lambda$-curve]\label{def:lambda curve}
For each pair of distinct points $x$, $y\in D$ and $\lambda>0$, we set $$c_{x,y}^\lambda:=\min\left\{1+\frac{\lambda}{2k_D(x,y)},\frac{9}{8}\right\}.$$ An arc $\gamma\subset D$ with end points $x$ and $y$ is called a $\lambda$-\textit{curve} if it is a $c_{x,y}^\lambda$-quasigeodesic in $D$. We use $\Lambda_{xy}^{\lambda}(D)$ to stand for the set of all $\lambda$-curves in $D$ with end points $x$ and $y$ and use $\gamma_{xy}^\lambda$ to denote a general $\lambda$-curve in $\Lambda_{xy}^{\lambda}(D)$. Throughout this paper, we always assume that $\lambda\in (0,\frac{1}{4}]$.
\edefe


By \cite[Theorem 3.3]{Vai4}, we know that for $x,y\in D$ and $\lambda>0$, if $x\not=y$, then $\Lambda_{xy}^{\lambda}(D)\not=\emptyset$.
For $c_0\geq 1$,
	we use $Q^{c_0}(D)$ to denote all $c_0$-quasigeodesics in $D$ and $Q^{c_0}_{xy}(D)$ for $c_0$-quasigeodesics in $D$ with end points $x$ and $y$. When the dependence on the domain $D$ is clear, we often drop it in the notations $S^{c,h}(D)$, $S^{c,h}_{xy}(D)$, $\Lambda_{xy}^\lambda(D)$, $Q^{c_0}(D)$ or $Q^{c_0}_{xy}(D)$. For instance, we often write $S^{c,h}$ and $S^{c,h}_{xy}$, instead of  $S^{c,h}(D)$ and $S^{c,h}_{xy}(D)$

	For $\lambda$-curves, we have the following two elementary properties.
	
	\blem[{\cite[Lemma 2.6]{Guo-H-W-2025}}]\label{lem-2.2}
	For any two distinct points $x\not=y\in D$, every $\lambda$-curve $\gamma\in \Lambda_{xy}^{\lambda}(D)$ is $\lambda$-short.
	\elem

	\blem[{\cite[Lemma 2.7]{Guo-H-W-2025}}]\label{9-27-1}
	Let $x\neq y$ be two distinct points in $D$ and fix a $\lambda$-curve $\gamma\in \Lambda_{xy}^\lambda(D)$. Then for any $z_1\not=z_2\in\gamma$, the subarc $\gamma[z_1,z_2]$ of $\gamma$ belongs to $ \Lambda_{z_1z_2}^{2\lambda}(D)$.
	\elem
	
	Next, we show that for a given special John curve (cone curve), we can construct an uniform curve with the same end points, that is also a quasigeodesic.
	\begin{lem}\label{thm3} For $u$ and $v$ in $D$, if there is an arc
		$\gamma$ in $D$ joining $u$ and $v$ such that for all $w\in\gamma$, $$\ell(\gamma[u,w])\leq c d_{D}(w),$$  then there exists a $c_1$-uniform arc $\beta\in Q^{c_1}_{uv}(D)$ with $c_1=c_1(c)$. Moreover, for all $w\in \beta$, \beq\label{eq-l-1}\ell(\beta[u,w])\leq c_1 d_D(w)\eeq
		and \be\label{eq-l-2-0}k_D(u,v)\leq 4c_1\log\Big(1+\frac{\ell(\beta)}{d_D(u)}\Big).\ee
	\end{lem}
	
	\bpf
	For $u$ and $v$ in $D$, let $\gamma\subset D$ be an arc joining $u$ and $v$ such that for all $w\in\gamma$, $$\ell(\gamma[u,w])\leq c d_{D}(w).$$ Then we infer from \cite[Lemma 2]{HLPW} (or \cite[Lemma 2.2]{Hl}) 
	that there exists a simply connected domain $D_{1}\subset D$ such that $u, v\in D_1$ and $D_1$ is a 
	$\varrho$-uniform domain with $\varrho=\varrho(c)$.
	
	Let $\beta\in Q_{uv}^3(D_1)$. Then it follows from  \cite[Theorem 6.19]{Vai4} that $\beta$ is a $\varrho_1$-uniform curve in $D_{1}$ with $\varrho_1=\varrho_1(\varrho)=\varrho_1(c)$, and so also a $\varrho_1$-uniform curve in $D$. Moreover, it follows from \cite[Theorem 2.42]{Vai7} that $\beta\in Q_{uv}^{\rho_2}(D)$ for some $\rho_2=\rho_2(\rho)=\rho_2(c)$.

	Let $w_0\in\beta$ bisect the length of $\beta$. 
	If $w\in \beta[u, w_0]$, then we have
	\[
	\ell(\beta[u,w])=\min\{\ell(\beta[u,w]),\ell(\beta[w,v])\}\leq \varrho_1d_D(w).
	\]
	If $w\in \beta[w_0,v]$, then we have 
	$$|v-w|\leq \ell(\beta[v,w])=\min\{\ell(\beta[u,w]),\ell(\beta[w,v])\}\leq \varrho_1 d_D(w)$$
	and so Lemma \ref{lem-3-1}(3) gives  
	\beqq
	d_D(v)\leq 2\varrho_1d_D(w).
	\eeqq
	This, together with the uniformity of $\beta$ and assumption of the lemma, implies that 
	$$\ell(\beta[u,w])\leq \ell(\beta)\leq \varrho_1|u-v|\leq \varrho_1\ell(\gamma)\leq c \varrho_1d_D(v)\leq 2c \varrho_1^2d_D(w).$$
	Thus \eqref{eq-l-1} holds if we choose $c_1=\max\{2c\varrho_1^2,\rho_2\}$. 
	
	The estimate \eqref{eq-l-2-0} follows then from \eqref{eq-l-1} and Lemma \ref{lem-3-1}$(2)$.
	
	\epf

	\subsection{Gromov hyperbolic and Rips spaces}
	
	There are many equivalent definitions of Gromov hyperbolicity and we shall mainly follow the one used in \cite{Vai10}. The following definition is equivalent to the one described in the introduction. 
	\begin{defn}\label{def:Gromov hyperbolicity}
		A domain $D\varsubsetneq E$ is called {\it $\delta$-Gromov hyperbolic}, $\delta>0$, if $(D,k_D)$ is $\delta$-Gromov hyperbolic.
		This means that for all $x, y, z, p\in D$, $$(x|z)_p\geq \min\{(x|y)_p,(y|z)_p\}-\delta,$$ where $(x|y)_p$ is the Gromov product defined by $$2(x|y)_p=k_D(p,x)+k_D(p,y)-k_D(x,y).$$
	\end{defn}
	
	It is well-known that every inner uniform domain is Gromov hyperbolic; see \cite[Theorem 1.11]{BHK}. 
	Since Gromov hyperbolicity is preserved by CQH mappings \cite[Theorem 3.18]{Vai10}, it follows that the quasiconformal image of an inner uniform domain in $\IR^n$ is Gromov hyperbolic. In particular, this implies by the Riemann mapping theorem that every simply connected proper subdomain of the complex plane is Gromov hyperbolic.

	Let $x_1,$ $x_2$ and $x_3$ be a triple of points in $D$. For each $i\in \{1,2,3\}$, let $\alpha_i$ denote an arc joining $x_{i}$ and $x_{i+1}$, where $x_4=x_1$.
	We write $\Delta=(\alpha_1,\alpha_2, \alpha_3)$ to represent a {\it triangle} in $D$ with sides $\alpha_1$, $\alpha_2$, $\alpha_3$ and the points $x_1,x_2,x_3$ are called the {\it vertices} of $\Delta$. A triangle $\Delta$ in $D$ is called {\it $h$-short} if all its three sides are $h$-short.

	%
	
	Next, we recall the definition of Rips space introduced in \cite[Section 2.26]{Vai10}. In what follows, we use the notation $k_D(w,\alpha)$ to denote the quasihyperbolic distance between the point $w\in D$ and the curve $\alpha$, i.e., $k_D(w,\alpha)=\inf\{k_D(w,v):\; v\in \alpha\}$.
	\bdefe[Rips space] \label{def12-1}
	Let $h\geq 0$ be a constant. If there is a constant $C>0$ such that for every $h$-short triangle $\Delta=(\alpha_1,\alpha_2, \alpha_3)$ in $D$ and for all $w\in \alpha_i$, $i\in \{1,2,3\}$, with $\alpha_{4}=\alpha_{1}$ and $\alpha_{5}=\alpha_{2}$, it holds
	$$k_D(w, \alpha_{i+1}\cup \alpha_{i+2})\leq C,$$
	then we call $(D, k_D)$ a {\it $(C,h)$-Rips space}.
	\edefe
	
	The following lemma gives uniform control on distance between ``close" points on short triangle in a Rips space.
	
	\begin{lem}[{\cite[Lemma 2.21]{Guo-H-W-2025}}]\label{lem03}  Suppose that $(D, k_D)$ is a $(C,h)$-Rips space. Then for every $h$-short triangle $\Delta=(\alpha_1,\alpha_2, \alpha_3)$ in $D$, and for each $i\in \{1,2,3\}$, there is $w_i\in \alpha_i$ such that for any $i\not=j\in \{1,2,3\}$,
		$$k_D(w_i,w_j)\leq 3C.$$ 
	\end{lem}

	\begin{lem}[{\cite[Theorem 2.35]{Vai10}}]\label{lem01} Let $D\subset E$ be a $\delta$-Gromov hyperbolic domain. Then $(D, k_D)$ is a $(C_0,h_0)$-Rips space for each $h_0> 0$ with $C_0$ depending only on $\delta$ and $h_0$.  
		In particular, $(D, k_D)$ is a $(C_0,h_0)$-Rips space for some $C_0$ and $h_0$ depending only on $\delta$. 
	\end{lem}
	
	The next lemma shows that the quasihyperbolic distance between a $\lambda$-curve and a quasigeodesic with the same end points is bounded, quantitatively.  
	\blem\label{Thm-23(1)}
	Let $D\subset E$ be $\delta$-Gromov hyperbolic and $x,y\in D$. Then there exists a constant $M_0=M_0(c,\delta)>0$ such that 
	\begin{itemize}
		\item for every $\gamma\in\Lambda_{xy}^{\lambda}(D)$ and $u\in \gamma$, $\beta\in Q_{xy}^{c}(D)$, $k_D(u, \beta)\leq M_0$;
		
		\item for every $\gamma\in\Lambda_{xy}^{\lambda}(D)$, $\beta\in Q_{xy}^{c}(D)$ and $v\in \beta$, $k_D(v,\gamma)\leq M_0$.
	\end{itemize}
	\elem
	\begin{proof}
		By \cite[Theorem 3.11]{Vai10}, we know that in a $\delta$-Gromov hyperbolic domain, if an $h$-short arc $\gamma$ and a $c$-quasigeodesic $\beta$ have the same end points, then the Hausdorff distance (with respect to the quasihyperbolic metric) between
		these two arcs is bounded by a constant which depends only on $c$, $h$ and $\delta$. Then the conclusions follow from the above fact by noticing Lemma \ref{lem-2.2}.
	\end{proof}

	The following result gives elementary estimate for ``quasigeodesic" triangles.
	\begin{lem}\label{lem24-1} 
		Let $D\subset E$ be $\delta$-Gromov hyperbolic, and $x,y,z\in D$. Then there exists a constant $A_0=A_0(\delta,c)>0$ such that for each $\gamma_1\in Q_{xy}^c(D)$, $\gamma_2\in Q_{xz}^c(D)$ and $\gamma_3\in Q_{yz}^c(D)$,  there are $w_i\in \gamma_i$, $i\in \{1,2,3\}$, such that for any $i\not=j\in \{1,2,3\}$,
		$$k_D(w_i,w_j)\leq A_0.$$ 
		Indeed, we can choose $A_0=3C_0+2M_0$, where $C_0=C_0(\delta)$ is the constant given by Lemma \ref{lem01} and $M_0=M_0(\delta,c)$ is the constant given by Lemma \ref{Thm-23(1)}.
	\end{lem}
	
	\bpf Fix $\beta_1\in \Lambda_{xy}^\lambda$, $\beta_2\in \Lambda_{yz}^\lambda$ and $\beta_3\in \Lambda_{xz}^\lambda$. Since $D$ is $\delta$-Gromov hyperbolic, Lemma \Ref{lem01} shows that $(D, k_D)$ is a $(C_0,h_0)$-Rips space for some $C_0=C_0(\delta)$ and $h_0=h_0(\delta)$, and thus it follows from Lemma \ref{lem03} that there exists $u_i\in \beta_{i}$ for each $i=1,2,3$, such that for each $i\not= j\in\{1,2,3\}$,
	\be\label{24-01-22-4}k_D(u_i,u_j)\leq 3C_0.\ee
	Further, by Lemma \ref{Thm-23(1)}, there are $w_1\in \gamma_{1}$,  $w_2\in \gamma_{2}$ and $w_3\in \gamma_{3}$
	such that for each $i\in\{1,2,3\}$,
	$$k_D(w_i,u_i)\leq M_0,$$
	which, together with (\ref{24-01-22-4}), shows that
	$$k_D(w_i,w_j)\leq k_D(w_i,u_i)+k_D(u_i,u_j)+k_D(w_j,u_j)\leq 3C_0+2M_0.$$
	\epf
	
	The next lemma compares the length of quasigeodesics with one common end point in Gromov hyperbolic spaces. 
	\begin{lem}[{\cite[Lemma 2.27]{Guo-H-W-2025}}]\label{main-lem-2} Suppose that $D$ is $\delta$-Gromov hyperbolic, $x_1$, $x_2$ and $x_3$ are points in $D$, and that $M_0$ is the constant given by Lemma \ref{Thm-23(1)}. Suppose further that there is a constant $\tau\geq 0$ such that $k_D(x_2,x_3)\leq\tau$. If $$\min\{k_D(x_1,x_2),k_D(x_1,x_3)\}> 1,$$
		then for any $\vartheta\in Q_{x_1x_2}^{c}(D)$ and $\gamma\in Q_{x_1x_3}^{c}(D)$, we have
		$$\max\left\{\ell(\vartheta),\ell(\gamma)\right\}\leq 4e^{2c^2(2\tau+3M_0)+\frac{4}{3}M_0}\min\left\{\ell(\vartheta),\ell(\gamma)\right\}.$$ 
	\end{lem}

	With Lemma \ref{main-lem-2} at hand, we are able to prove the following technical result, which shall be important in the proof of Theorem \ref{thm-1}.
	\begin{lem}\label{lem02}
		Let $D$ be $\delta$-Gromov hyperbolic and $x$, $y$ and $z$ be three points in $D$. Fix $\gamma_1\in Q_{xy}^c(D)$, $\gamma_2\in Q_{xz}^c(D)$ and assume $k_D(y,z)\leq \tau$ with $\tau\geq 0$. Then we have
		\begin{enumerate}
			\item\label{lem02-1}
			For each $w\in \gamma_1$, there exists $v\in \gamma_2$ such that
			$$k_D(w,v)\leq \frac{3}{4}+C_0+2M_0+\tau.$$
			
			\item\label{lem02-2} If for each $w\in \gamma_{1}$, it holds $\ell(\gamma_{1}[x,w])\leq cd_D(w)$, then for each $u\in \gamma_{2}$, we have 
			$\ell(\gamma_{2}[x,u])\leq \theta d_D(u)$ with $\theta=e^{4c^2(2+2C_0+4M_0+\tau)}c$. 
		\end{enumerate}
		Here $C_0$ and $M_0$ are two constants given by Lemmas \ref{lem01} and \ref{Thm-23(1)}, respectively.
	\end{lem}
	
	\bpf
	(1) Fix $\alpha_1\in \Lambda_{xy}^\lambda$, $\alpha_2\in \Lambda_{xz}^\lambda$ and $\alpha_3\in \Lambda_{yz}^\lambda$. It follows from Lemma \ref{Thm-23(1)} that there exists $w_1\in\alpha_1$ such that
	\be\label{24-01-22-1}k_D(w,w_1)\leq M_0.\ee
	Since $D$ is $\delta$-Gromov hyperbolic, Lemma \Ref{lem01} implies that $(D,k_D)$ is a $(C_0,h_0)$-Rips space. It follows 
	$$k_D(w_1,\alpha_2\cup \alpha_3)\leq C_0.$$
	Select $w_2\in \alpha_2\cup \alpha_3$ so that $$k_D(w_1,w_2)\leq \frac{1}{2}+C_0.$$
	
	If $w_2\in \alpha_3$, then take $v=z$, and we obtain from Lemma \ref{lem-2.2} and \eqref{24-01-22-1} that
	\[
	\begin{aligned}
		k_D(w,v)&\leq k_D(w,w_1)+k_D(w_1,w_2)+k_D(w_2,v)\\ &\leq \frac{1}{2}+C_0+M_0+\ell_k(\alpha_3)
		\leq \frac{3}{4}+C_0+M_0+k_D(y,z)\\&\leq \frac{3}{4}+C_0+M_0+\tau.
	\end{aligned}
	\]

	If $w_2\in \alpha_2$, then by Lemma \ref{Thm-23(1)}, there exists $v\in \gamma_2$ such that
	$$k_D(w_2,v)\leq M_0.$$
	This, together with \eqref{24-01-22-1}, gives
	$$k_D(w,v)\leq k_D(w,w_1)+k_D(w_1,w_2)+k_D(w_2,v)\leq \frac{1}{2}+C_0+2M_0.$$
	
	(2) Fix $u\in \gamma_2$. By statement \eqref{lem02-1}, there exists $w\in \gamma_1$
	such that
	\be\label{24-01-22-2}k_D(u,w)\leq\frac{3}{4}+C_0+2M_0+\tau.\ee
	This, together with \eqref{(2.2)}, implies that $$\log \frac{d_D(w)}{d_D(u)}\leq k_D(u,w)\leq\frac{3}{4}+C_0+2M_0+\tau,$$
	and thus we have 
	\be\label{24-01-22-3}d_D(w)\leq e^{\frac{3}{4}+C_0+2M_0+\tau}d_D(u).\ee
	
	If $\min\{k_D(x,u),k_D(x,w)\}\leq 1$, then by \eqref{24-01-22-2}, we have
	\beqq
	k_D(x,u)\leq \min\{k_D(x,u),k_D(x,w)\}+k_D(u,w)\leq\frac{7}{4}+C_0+2M_0+\tau.
	\eeqq
	Combing this with \eqref{(2.1)} gives
	$$\log\Big(1+\frac{\ell(\gamma_{2}[x,u])}{d_D(u)}\Big)\leq \ell_k(\gamma_{2}[x,u])\leq ck_D(x,u)\leq c\left(\frac{7}{4}+C_0+2M_0+\tau\right),$$
	from which we get
	$$\ell(\gamma_{2}[x,u])< e^{c(\frac{7}{4}+C_0+2M_0+\tau)}d_D(u).$$
	
	If $\min\{k_D(x,u),k_D(x,w)\}> 1$, then it follows from Lemma \ref{main-lem-2} (applied with $x=x_1$, $x_2=u$ and $x_3=w$) and statement \eqref{lem02-1} of the lemma that
	$$
	\begin{aligned}
		\ell(\gamma_{2}[x,u])&\leq e^{4c^2(\frac{3}{4}+C_0+3M_0+\tau)}\ell(\gamma_{1}[x,w])\leq e^{4c^2(\frac{3}{4}+C_0+3M_0+\tau)}cd_D(w)\\
		&\stackrel{\eqref{24-01-22-3}}{\leq}e^{4c^2(1+2C_0+4M_0+\tau)}cd_D(u)\leq \theta d_D(u).
	\end{aligned}
	$$
	\epf

	\section{Dimension-free cone property of quasigeodesics} \label{sec-2}
	In this section, we establish quantitative cone property of quasigeodesics in a Gromov hyperbolic John domain. 
	\begin{thm}\label{thm-1} 
		Let $E$ be a Banach space and let $D\subsetneq E$ be a $\delta$-Gromov hyperbolic $c$-John domain. Then every $c_0$-quasigeodesic $\gamma\subset D$ is a $c_2$-cone arc with $c_2=c_2(\delta,c,c_0)$.
	\end{thm}

	\begin{proof}
		We shall show that every $c_0$-quasigeodesic in a $\delta$-Gromov hyperbolic $c$-John domain is a $c_2$-cone arc with $c_2=e^{20c_0^2c_1^2A_1}$, where $c_1=c_1(c)>1$ is the constant from Lemma \ref{thm3} and $A_1=A_1(\delta,c_0)$ (cf. Lemma \ref{lem24-1}).
		
		Since $D$ is $\delta$-Gromov hyperbolic, Lemma \Ref{lem01} yields that for some constants $C_0=C_0(\delta)$ and $h_0=h_0(\delta)$, $(D, k_D)$ is a $(C_0,h_0)$-Rips space. 
		Fix $z_1,z_2\in D$, and let $\gamma\in Q_{z_1z_2}^{c_0}(D)$. 
		Since $D$ is a $c$-John domain, there is a $c$-cone arc $\alpha$ in $D$  joining $z_1$ and $z_2$. Let $z_0\in \alpha$ be the point bisecting the length of $\alpha$. Then Lemma \ref{thm3} implies that there exist two $c_1$-uniform arcs $\alpha_1\in Q_{z_0z_1}^{c_1}(D)$ and $\alpha_2\in Q_{z_0z_2}^{c_1}(D)$ with $c_1=c_1(c)$. Moreover, for each $i\in\{1,2\}$ and all $w\in \alpha_i$, it holds 
		\beq\label{eq-ne-1}\ell(\alpha_i[z_i,w])\leq c_1 d_D(w).\eeq
		
		Now, we claim that for $c_2=e^{20c_0^2c_1^2A_1}$ given above, it holds that for all $w\in \gamma$,
		\beq\label{eq-l-0}\min\{\ell(\gamma[z_1,w]),\ell(\gamma[z_2,w])\}\leq c_2 d_D(w).\eeq
		
		Set $\alpha_3=\gamma$. Then all the $\alpha_i$ are $(c_0+c_1)$-quasigeodesics in $D$. It follows from Lemma \ref{lem24-1} that for each $i\in\{1,2,3\}$, there exists $x_i\in\alpha_i$, such that
		for each $s\not=t\in\{1,2,3\}$,
		\be\label{24-01-23-1}k_D(x_s,x_t)\leq A_1,\ee
		where $A_1=A_1(\delta, c_0+c_1)$.
		
		Next we prove that the two subcurves $\alpha_3[z_1,x_3]$ and $\alpha_3[z_2,x_3]$ satisfy \eqref{eq-l-0}.
		
		For each $x\in\alpha_3[z_1,x_3]$, it follows from Lemma \ref{lem02}(\ref{lem02-2}), \eqref{eq-ne-1} and \eqref{24-01-23-1} that
		$$\ell(\alpha_3[z_1,x])< e^{20c_0^2c_1^2A_0}d_D(x)=c_2d_D(x).$$
		For each $x\in\alpha_3[z_2,x_3]$, it follows again from Lemma \ref{lem02}(\ref{lem02-2}), \eqref{eq-ne-1} and \eqref{24-01-23-1} that
		$$\ell(\alpha_3[z_2,x])< e^{20c_0^2c_1^2A_0} d_D(x)=c_2d_D(x).$$
		This proves our claim \eqref{eq-l-0} and thus completes the proof of Theorem \ref{thm-1}.
	\end{proof}
	
	\begin{rem}\label{rmk:on Theorem 3}
		We would like to point out that Theorem \ref{thm-1} improves on \cite[Proposition 7.12]{BHK} and \cite[Remark 3.10]{Guo15} as we do not need the ``boundedness"  assumption and our estimate is dimension free. 
	\end{rem}
	
	As a corollary of Theorem \ref{thm-1}, we obtain a generalization of \cite[Theorem 1]{Li} 
	via a new approach.

	\begin{cor}\label{coro2}
		Suppose that $E$ is a Banach space and $D\subset E$ is a $c$-John domain, which
		homeomorphic to a $\delta$-Gromov hyperbolic domain via a $(M,C)$-CQH mapping. Then  every $c_0$-quasigeodesic in $D$ is a $b$-cone arc with $b=b(\delta,c,c_0,M,C)$.
	\end{cor}

	\section{Proof of Theorem \ref{thm:main} } \label{sec-3}	
	\subsection{Preparation for the proof}
	Let $E$ be a Banach space and $D\subsetneq E$ a $\delta$-Gromov hyperbolic $c$-John domain. By Theorem \ref{thm-1}, a $c_0$-quasigeodesic $\gamma$ in $D$ satisfies condition (1) of Definition \ref{def:uniform domain}. To prove Theorem \ref{thm:main}, it suffices to prove that this quasigeodesic  $\gamma$ satisfies the Gehring-Hayman inequality \eqref{eq:def for inner quasiconvexity}. 
	
	In the following discussion, we fix $x,y\in D$. Let $\gamma$ be a $c_0$-quasigeodesic joining $x$ and $y$ in $D$ and let $\alpha$ be an arbitrary curve joining $x$ and $y$ in $D$. As in \cite{Guo-H-W-2025}, we shall prove Theorem \ref{thm:main} by a contradiction argument.
	
	Suppose on the contrary that for some large $\kappa\gg 1$, it holds
	\be\label{lem2-eq-2}
	\ell(\gamma)\geq \kappa\ell(\alpha).
	\ee
	Under this assumption, we have $k_D(x,y)> 1$. Indeed, if $k_D(x,y)\leq 1$,  then we obtain from \cite[Lemma 2.31]{Guo-H-W-2025}  that $\gamma$ is a  $18c_0e^{9c_0}$-uniform arc, which implies
	$$\ell(\gamma)\leq 18c_0e^{9c_0}|x-y|\leq 18c_0e^{9c_0}\ell(\alpha).$$
	This contradicts with \eqref{lem2-eq-2} if we select $\kappa>18c_0e^{9c_0}$. Furthermore, for each $w\in\alpha$, we know from \cite[Lemma 2.2]{Vai7} that
	\be\label{lem2-eq-3}d_D(w)\leq 2\ell(\alpha).\ee
	
	Let $y_0\in\gamma$ bisect the length of $\gamma$. Then Theorem \ref{thm-1} implies that for all $w\in \gamma$,
	$$	\min\{\ell(\gamma[x,w]),\ell(\gamma[y,w])\}\leq c_2 d_D(w),$$
	and thus by \eqref{lem2-eq-2}, we have
	\be\label{24-01-28-2} 
	d_D(y_0)\geq (2c_2)^{-1}\ell(\gamma)> \frac{\kappa-2}{2c_2}\ell(\alpha).
	\ee

	\subsection{$\zeta$-sequence on $\alpha$}\label{sec-3.1}

	\begin{figure}[htbp]
		\begin{center}
			\begin{tikzpicture}[scale=2.5]
				\draw (-2,0) to [out=60,in=180] (-1,0.8) to [out=0,in=120] (0,0) to [out=300,in=180] (1.7,-0.7);
				\draw (1.7,-0.7) to [out=0,in=-126] (3.304,0.1);
				\filldraw  (-2,0) circle (0.02);
				\node at(-1.85,-0.15) {\small $x\!\!=\!x_1\!\!=\!\zeta_1$};
				\filldraw  (-1.4,0.7) circle (0.02);
				\node[rotate=-60, xshift=-0.4cm,  yshift=-0.1cm] at(-1.18,0.3) {\small $x_2=\zeta_2$};
				\filldraw  (-0.7,0.745)node[below, rotate=-15] {\small $x_3=\zeta_3$} circle (0.02);
				\filldraw  (-0.25,0.392)node[below left, xshift=0.2cm] {\small $x_{i-1}=\zeta_{i-1}$} circle (0.02);
				\filldraw  (0.2,-0.265)node[below, xshift=0.3cm] {\small $x_{i}\!\!=\!\zeta_{i}\;\;\;\;\;\;\;\;\;\;\;\;\;\;\;\;$} circle (0.02);
				\filldraw  (0.85,-0.622)node[below, rotate=15] {\small $x_{i+1}\!\!=\zeta_{i+1}\;\;\;\;\;\;\;\;$} circle (0.02);
				\filldraw  (1.45,-0.695)node[below, rotate=0] {\small $\;\;\;\;\;\;\;\;\;\;\;\;x_{m-1}\!\!=\zeta_{m-1}$} circle (0.02);
				\filldraw  (2.5,-0.56)node[below] {\small $\;\;\;x_m$} circle (0.02);
				\filldraw  (3.304,0.1)node[below, rotate=50] {\small $y\!=\!x_{m\!+\!1}\!\!=\!\zeta_m$} circle (0.02);
				
				\draw (-2,0) to [out=95,in=180] (-1.7,1) to [out=0,in=110] (-1.4,0.7);
				\filldraw  (-1.95,0.865)node[below right] {\tiny $w_1$} circle (0.02);
				\node[] at(-1.82,1.1) {\tiny $\alpha_1\in \Lambda_{x_1x_2}^{\lambda}$};
				
				\draw (-1.4,0.7) to [out=85,in=180] (-1.05,1.3) to [out=0,in=100] (-0.7,0.745);
				\filldraw  (-1.25,1.213)node[below right] {\tiny $w_2$} circle (0.02);
				\node[above] at(-1.05,1.3) {\tiny $\alpha_2\in \Lambda_{x_2x_3}^{\lambda}$};
				
				\draw (-0.25,0.392) to [out=85,in=180] (-0.1,0.6) to [out=0,in=90] (0.2,-0.265);
				\filldraw  (0.15,0.43)node[below left, xshift=0.125cm] {\tiny $w_{i-1}$} circle (0.02);
				\node[above, rotate=-10] at(0.037,0.55) {\tiny $\alpha_{i-1}\in \Lambda_{x_{i-1}x_i}^{\lambda}$};

				\draw (0.2,-0.265) to [out=80,in=185] (0.5,0.1) to [out=0,in=90] (0.85,-0.622);
				\filldraw  (0.65,0.06)node[below] {\tiny $w_{i}$} circle (0.02);
				\node[above, rotate=-10, xshift=0.225cm, yshift=-0.06cm] at(0.5,0.1) {\tiny $\alpha_i\!\!\in \!\! \Lambda_{x_ix_{i+1}}^{\lambda}$};

				\draw (1.45,-0.695) to [out=67,in=180] (2,-0.35) to [out=0,in=135] (2.5,-0.555);
				\filldraw[xshift=0.6cm,yshift=-0.085cm]  (1.5,-0.27)node[below] {$\eta_1$} circle (0.02);
				\node[rotate=5, xshift=1.5cm ,yshift=-0.4cm] at(1.45,-0.135) {\tiny $\alpha_{m-1}\!\!\in \!\! \Lambda_{x_{m-1}x_m}^{\lambda}$};
				
				\draw (2.5,-0.555) to [out=95,in=220] (2.72,0) to [out=40,in=170] (3.304,0.1);
				\filldraw[xshift=0.6cm]  (2.12,0) circle (0.02);
				\node[xshift=01.5cm] at(2.2,-0.1) {\small $\eta_2$};
				\node[xshift=1.5cm, rotate=30] at(2.1,0.15)  {\tiny $\alpha_{m}\!\!\in \!\! \Lambda_{x_{m}x_{m+1}}^{\lambda}$};
				
				\draw (1.45,-0.695) to [out=80,in=200] (2.2,0.4) to [out=20,in=135] (3.304,0.1);
				\filldraw[xshift=0.6cm,yshift=-0.12cm]  (1.2,0.257)node[below, xshift=0.3cm,] {\tiny $w_{m-1}$} circle (0.02);
				\filldraw[xshift=0.6cm,yshift=-0.12cm]  (1.55,0.5)node[below] {$\eta_3$} circle (0.02);
				\node[xshift=1.3cm,yshift=-0.15cm, rotate=-5] at(2.1,0.65) {\tiny $\varsigma_m=\alpha_{m-1}\!\in \! \Lambda_{x_{m-1}x_{m+1}}^{\lambda}$};
				
				\node at(0.38,-0.6) {\large $\alpha$};
				
				
				\filldraw  (0.5,4.5)node[above] { $y_0$} circle (0.02);
				\coordinate (Y) at (0.5,4.5);
				
				\node at(1.7,4.4) { $\varsigma\in \Lambda_{xy_1}^{\lambda}$};
				
				\draw (-2,0) to [out=110,in=175] (Y) to [out=-5,in=90] (3.304,0.1);
				\filldraw  (-1.94,2.6)node[below right] {\tiny $w_{1,1}$} circle (0.02);
				\node[left] at(-1.65,3.2) {\small $\gamma_1$};
				\filldraw  (-1.31,3.7)node[below right] {\tiny $w_{0,1}$} circle (0.02);
				\filldraw[color=red]  (-1,4.004)node[below] {\tiny $\xi_1$} circle (0.02);
				
				\draw (-1.4,0.7) to [out=100,in=200] (Y);
				\filldraw  (-1.445,1.89)node[right] {\tiny $w_{1,2}$} circle (0.02);
				\filldraw  (-1.355,2.4)node[right] {\tiny $w_{2,1}$} circle (0.02);
				\node[left] at(-1.13,3.02) {\small $\gamma_2$};
				\filldraw  (-0.825,3.53)node[below right] {\tiny $w_{0,2}$} circle (0.02);
				\filldraw[color=red]  (-0.55,3.85)node[below] {\tiny $\xi_2$} circle (0.02);
				
				\draw (-0.7,0.745) to [out=90,in=220] (Y);
				\filldraw  (-0.687,1.695)node[right] {\tiny $w_{2,2}$} circle (0.02);
				\filldraw  (-0.65,2.21)node[right] {\tiny $w_{3,1}$} circle (0.02);
				\node[left] at(-0.54,2.87) {\small $\gamma_3$};
				\filldraw  (-0.342,3.41)node[below] {\tiny $\;\;\;\;\;\;w_{0,3}$} circle (0.02);
				\filldraw[color=red]  (-0.15,3.77)node[below] {\tiny $\;\;\xi_3$} circle (0.02);
				
				\draw (-0.25,0.392) to [out=88,in=248] (Y);
				\filldraw  (-0.188,1.61)node[right] {\tiny $w_{i\!-\!2,2}$} circle (0.02);
				\filldraw  (-0.135,2.13)node[right] {\tiny $w_{i-1,1}$} circle (0.02);
				\node[left] at(-0.04,2.79) {\small $\gamma_{i-1}$};
				\filldraw  (0.105,3.35)node[below] {\tiny $\;\;\;\;\;\;\;\;\;\;\;w_{0,i-1}$} circle (0.02);
				\filldraw[color=red]  (0.21,3.703)node[below, xshift=0.1cm] {\tiny $\;\;\;\xi_{i-1}$} circle (0.02);
				
				\draw (0.2,-0.265) to [out=80,in=-90] (Y);
				\filldraw  (0.44,1.58)node[right] {\tiny $w_{i-1,2}$} circle (0.02);
				\filldraw  (0.47,2.1)node[right] {\tiny $w_{i,1}$} circle (0.02);
				\node[left] at(0.54,2.77) {\small $\gamma_{i}$};
				\filldraw  (0.5,3.33)node[below] {\tiny $\;\;\;\;\;\;\;\;\;w_{0,i}$} circle (0.02);
				\filldraw[color=red]  (0.5,3.7)node[below] {\tiny $\;\;\;\;\;\xi_{i}$} circle (0.02);
				
				\draw (0.85,-0.622) to [out=85,in=-70]coordinate[pos=0.41] (wi2) coordinate[pos=0.52] (wi+11) node [left, pos=0.656, xshift=0.1cm]{\small $\gamma_{i\!+\!1}$} coordinate[pos=0.78] (w0i+1) coordinate[pos=0.85] (xii+1)  (Y);
				\filldraw  (wi2)node[right, xshift=-0.25cm] {\tiny $\;\;\;w_{i,2}$} circle (0.02);
				\filldraw  (wi+11)node[right, xshift=-0.25cm] {\tiny $\;\;\;w_{i\!+\!1,1}$} circle (0.02);
				\filldraw  (w0i+1)node[below right, xshift=-0.3cm] {\tiny $\;\;\;w_{0,i\!+\!1}$} circle (0.02);
				\filldraw[red]  (xii+1)node[below right , xshift=-0.25cm] {\tiny $\;\;\;\xi_{i\!+\!1}$} circle (0.02);
				
				\draw (1.45,-0.695) to [out=85,in=-45] (Y);
				\filldraw[xshift=0.652cm]  (0.965,1.605)node[right] {\tiny $w_{m\!-\!2,2}$} circle (0.02);
				\filldraw[xshift=0.625cm]  (0.965,2.13)node[right] {\tiny $w_{m\!-\!1\!,1}$} circle (0.02);
				\node[left, xshift=1.43cm] at(0.98,2.75) {\small $\gamma_{m\!-\!1}$};
				\filldraw[xshift=0.465cm]  (0.825,3.36)node[below right] {\tiny $w_{0,m\!-\!1}$} circle (0.02);
				\filldraw[color=red]  (1.115,3.708)node[below right, xshift=0.08cm, yshift=0.1cm] {\tiny $\xi_{m\!-\!1}$} circle (0.02);
				
				\filldraw[xshift=0.439cm]  (2.67,1.8)node[left] {\tiny $w_{m\!-\!1,2}$} circle (0.02);
				\node[left,xshift=0.8cm] at(2.4,2.9) {\small $\gamma_{m}$};
				\filldraw[xshift=0.23cm]  (2.1,3.5)node[below left] {\tiny $w_{0,m}$} circle (0.02);
				\filldraw[color=red]  (1.96,3.89)node[below] {\tiny $\xi_m\;\;$} circle (0.02);
				
				\draw[color=red] (-1,4.004) to [out=30,in=110]node [above, pos=0.5, xshift=0.08cm, yshift=-0.1cm]{\tiny $\beta_1$} (-0.55,3.85);
				\draw[color=red] (-0.55,3.85) to [out=30,in=110]node [above, pos=0.5, xshift=0.08cm, yshift=-0.1cm]{\tiny $\beta_2$} (-0.15,3.77);
				\draw[dashed, color=red] (-0.15,3.77) to [out=30,in=110] (0.21,3.703);
				\draw[color=red] (0.21,3.703) to [out=45,in=118]node [above, pos=0.5, xshift=0.06cm, yshift=-0.1cm]{\tiny $\beta_{i-1}$} (0.5,3.7);
				\draw[color=red] (0.5,3.7) to [out=45,in=118]node [above, pos=0.5, xshift=0cm, yshift=-0.1cm]{\tiny $\beta_i$} (xii+1);
				\draw[dashed, color=red] (xii+1) to [out=45,in=118] (1.115,3.708);
				\draw[color=red] (1.115,3.708) to [out=60,in=145] (1.96,3.89);
				
				\node[color=red] at(1.35,4.05) {\tiny $\beta_0$};
			\end{tikzpicture}
		\end{center}
		\caption{Illustration for the construction used in proof of Theorem \ref{thm:main}} \label{fig3}
	\end{figure}
	
	In this section, we shall recall from \cite[Section 5.2]{Guo-H-W-2025} the construction of a finite sequence of points on $\alpha$, called $\zeta$-sequence, with good controls on the length of $\lambda$-curves between successive points. 
	
	Let $x_1=x$ and $\kappa_1=\kappa^{\frac{1}{3}}$. Then we have the following sequence of finitely many points on $\alpha$. 
	
	\begin{prop}\label{24-5-7-3}
		There exist an integer $m\geq 2$ and a sequence $\{x_i\}_{i=1}^{m+1}$ of points on $\alpha$ such that the following conclusions hold.
		
		$(a)$ For each $i\in \{2,\cdots,m\}$ and for any $\lambda$-curve $\alpha_{i-1}\in \Lambda_{x_{i-1}x_i}^{\lambda}(D)$, we have
		\beqq\label{lem2-eq-5}\kappa_1^{\frac{1}{4}}\ell(\alpha)\leq \ell(\alpha_{i-1})\leq \kappa_1^{\frac{1}{2}}\ell(\alpha).\eeqq

		$(b)$ For any $z\in \alpha[x_i,y]$, and for each $\lambda$-curve $\varsigma_{i-1}\in \Lambda_{x_{i-1}z}^{\lambda}(D)$,
		\beqq\label{1-16-1}
		\ell(\varsigma_{i-1})>\kappa_1^{\frac{1}{4}}\ell(\alpha).
		\eeqq
		
		$(c)$ There is an element  $\alpha_m\in \Lambda_{x_mx_{m+1}}^\lambda(D)$ such that
		\beqq\label{1-16-01}
		0<\ell(\alpha_m)\leq \kappa_1^{\frac{1}{2}}\ell(\alpha).
		\eeqq
		Here $x_1=x$ and $x_{m+1}=y$.
	\end{prop}
	\bpf  
	It follows from the assumption $k_D(x,y)>1$, Lemma \ref{main-lem-2} and (\ref{lem2-eq-2}) that for each $\gamma_{xy}\in \Lambda_{xy}^{\lambda}(D)$, 
	$$\ell(\gamma_{xy})>c_2^{-1}\ell(\gamma)\geq \frac{\kappa}{c_2}\ell(\alpha),$$
	which, together with (\ref{lem2-eq-3}), shows that
	the assumptions in \cite[Lemma 2.29]{Guo-H-W-2025} are satisfied for any $\lambda$-curve $\varsigma\in \Lambda_{xy}^{\lambda}(D)$ and $\alpha$. Then by
	\cite[Lemma 2.29]{Guo-H-W-2025}, the rest proof of Proposition \ref{24-5-7-3} follows from \cite[Proposition 5.4]{Guo-H-W-2025}.
	\epf
	\medskip

	\blem[{\cite[Lemma 5.5]{Guo-H-W-2025}}]\label{lem2-eq-6} There is an element in $\varsigma_m\in\Lambda_{x_{m-1}x_{m+1}}^\lambda$ such that
	$$\kappa_1^{\frac{1}{4}}\ell(\alpha)<\ell(\varsigma_m)\leq \kappa_1^{\frac{3}{4}}\ell(\alpha).$$
	\elem
	
	For convenience, we make the following notational conventions. \begin{enumerate}
		\item\label{24-5-7-1}
		for each $i\in\{1,2,\cdots,m-1\}$, let $\zeta_i=x_i$ and $\zeta_m=x_{m+1}$.
		\item\label{24-5-7-2}
		For each $i\in\{1,2,\cdots,m-2\}$, we fix an element $\alpha_i\in\Lambda_{\zeta_i\zeta_{i+1}}^\lambda$. Moreover, let $\alpha_{m-1}=\varsigma_m\in \Lambda_{x_{m-1}x_{m+1}}^{\lambda}$ be given by Lemma \ref{lem2-eq-6}, and $\alpha_m\in \Lambda_{x_mx_{m+1}}^{\lambda}$ be given by Proposition \ref{24-5-7-3}$(c)$.
	\end{enumerate}
	Then we see from Proposition \ref{24-5-7-3}$(a)$ and Lemma \ref{lem2-eq-6} that for each $i\in\{1,2,\cdots,m-2\}$,
	$$\kappa_1^{\frac{1}{4}}\ell(\alpha)<\ell(\alpha_{i})\leq \kappa_1^{\frac{1}{2}}\ell(\alpha)\;\;\mbox{and}\;\;\kappa_1^{\frac{1}{4}}\ell(\alpha)<\ell(\alpha_{m-1})\leq \kappa_1^{\frac{3}{4}}\ell(\alpha),$$
	and from Proposition \ref{24-5-7-3}$(b)$ for each $\zeta\in\alpha[\zeta_{i+1},y]$, and for each $\varsigma_i\in 
	\Lambda_{x_i\zeta}^\lambda$, 
	$$\ell(\varsigma_i)>\kappa_1^{\frac{1}{4}}\ell(\alpha).$$
	Furthermore, we infer from \eqref{(2.1)} and \eqref{lem2-eq-3} that for each $i\in\{1,2,\cdots,m-1\}$,
	$$\ell_{k}(\alpha_i)\geq \log\Big(1+\frac{\ell(\alpha_i)}{d_D(\zeta_i)}\Big)\geq \log\Big(1+\frac{\ell(\alpha_i)}{2\ell(\alpha)}\Big)> \frac{1}{4}\log \frac{\kappa_1}{16}.$$

	\subsection{$\lambda$-curves associated to $\zeta$-sequence} 
	The presentation of this section follows closely the one in \cite[Section 5.3]{Guo-H-W-2025}. 
	Fix any $\lambda$-curve $\varsigma\in \Lambda_{xy}^\lambda$ and write  $\gamma_1=\varsigma[\zeta_1,y_0]$ and $\gamma_{m}=\varsigma[y,y_0]$. For $i\in \{2,3,\cdots,$ $m-1\}$, let
	$\gamma_i\in \Lambda_{\zeta_iy_0}^{\lambda}$ be an element; see Figure \ref{fig3}.
	
	Since $(D,k_D)$ is a $(C_0,h_0)$-Rips space, it follows from Lemmas \ref{lem-2.2} and \ref{lem03} that for each $i\in \{1,2,\cdots,m-1\}$, there exist three points $w_i\in \alpha_i$, $w_{i,1}\in \gamma_i$, $w_{i,2}\in \gamma_{i+1}$ such that
	$$k_D(w_i,w_{i,1})\leq 3C_0,\;\;k_D(w_{i,1},w_{i,2})\leq 3C_0\;\;{\rm and}\;\; k_D(w_i,w_{i,2})\leq 3C_0.$$
	Furthermore, we have the following estimate for $\ell(\gamma_i[\zeta_i,w_{i,1}])$ and $\ell(\gamma_{i+1}[\zeta_{i+1},w_{i,2}])$. 
	\blem[{\cite[Lemma 5.6]{Guo-H-W-2025}}]\label{lemma:2022-10-24-4}
	For each $i\in \{1,2,\cdots,m-1\}$, we have
	\be\label{2022-10-24-4}
	\max\{\ell(\gamma_i[\zeta_i,w_{i,1}]),\ell(\gamma_{i+1}[\zeta_{i+1},w_{i,2}])\}\leq 4\kappa_1^{\frac{3}{4}}e^{3(5C_0+3M_0)}\ell(\alpha).
	\ee
	\elem
	
	Next, we are going to find one more point on every $\gamma_i$. Note that 
	$$\ell(\gamma_i)\geq d_D(y_0)-d_D(\zeta_i)\stackrel{\eqref{lem2-eq-3}+\eqref{24-01-28-2}}{\geq}  \Big(\frac{\kappa-2}{2c_2}-2\Big)\ell(\alpha)>2\kappa_1\ell(\alpha).$$
	Thus for each $i\in\{1,\cdots,m\}$, there is $w_{0,i}\in\gamma_i$ such that
	$$\ell(\gamma_i[\zeta_i,w_{0,i}])=\kappa_1\ell(\alpha).$$
	Furthermore, we infer from \eqref{2022-10-24-4} that $w_{0,1}\in\gamma_1[w_{1,1},y_0]$,
	and for each $i\in\{2,\cdots,m\}$,
	$$w_{0,i}\;\in \; \gamma_i[w_{i,1},y_0]\cap \gamma_i[w_{i-1,2},y_0]. $$  
	
	We have the following two basic estimates from \cite{Guo-H-W-2025}. 
	\blem[{\cite[Lemma 5.7]{Guo-H-W-2025}}]\label{lemma:2022-09-30-2} For each $i\in\{1,\cdots,m-1\}$,
	$$\min\left\{k_D(w_{i,1},w_{0,i}),\;\;k_D(w_{i,2},w_{0,{i+1}})\right\} \geq \frac{1}{5}\log \kappa_1.
	$$
	\elem
	
	\blem[{\cite[Lemma 5.8]{Guo-H-W-2025}}]\label{lemma:23-07-23-4} For each $i\in\{1,\cdots,m-1\}$,
	$$\min\{k_D(y_0,w_{i,1})-k_D(w_{0,i},y_0),\; k_D(y_0,w_{i,2})-k_D(w_{0,i+1},y_0)\}  \geq \frac{1}{5}\log \kappa_1-\frac{1}{2}.$$
	\elem

	Let $t_1\in\{1,\cdots,m\}$ be such that
	\be\label{11-12-1} k_D(y_0,w_{0,t_1})=\min_{i\in\{1,\cdots,m\}}\{k_D(y_0,w_{0,i})\}.\ee
	Then we have the following assertion.
	\blem[{\cite[Lemma 5.9]{Guo-H-W-2025}}]\label{cl-23-07-23}
	Suppose that there exist $i\in\{1,\cdots,m\}$ and $u_i\in \gamma_i[w_{0,i},y_0]$ such that  $k_D(y_0,u_i)=k_D(y_0,w_{0,t_1})$.
	\begin{enumerate}
		\item\label{23-07-23-6} If $2\leq i\leq m$, then there exists $u_{i-1}\in \gamma_{i-1}$ such that
		$$k_D(y_0,u_{i-1})=k_D(y_0,w_{0,t_1})\;\;\mbox{and}\;\; k_D(u_i,u_{i-1})\leq 2(1+5C_0+2M_0).$$
		
		\item\label{23-07-23-7} If $1\leq i\leq m-1$, then there exists $u_{i+1}\in \gamma_{i+1}$ such that
		$$k_D(y_0,u_{i+1})=k_D(y_0,w_{0,t_1})\;\;\mbox{and}\;\;  k_D(u_i,u_{i+1})\leq 2(1+5C_0+2M_0).$$
	\end{enumerate}
	\elem
	
	\subsection{$\xi$-sequence on the associated $\lambda$-curves}\label{subsec:xi sequence} 
	
	In this section, following \cite[Section 5.4]{Guo-H-W-2025}, we shall construct a new sequence of points. By the choice of $w_{0,t_1}$ in \eqref{11-12-1} and Lemma \ref{cl-23-07-23}, we know that
	for each $i\in\{1,\cdots,m\}$, there exists $\xi_i\in\gamma_i[w_{0,i},y_0]$ (see Figure \ref{fig3}) such that
	\be\label{2022-09-30-3}
	k_D(y_0,\xi_i)=k_D(y_0,w_{0,t_1})\;\;\mbox{and}\;\;k_D(\xi_i,\xi_{i+1})<2(1+5C_0+2M_0)\ee

	\blem[{\cite[Lemma 5.10]{Guo-H-W-2025}}]\label{cl3.2-1} $(i)$
	For each $i\in\{1,\cdots,m\}$ and any $z\in\gamma_i[\zeta_i,\xi_i]$, we get
	$$k_D(z,y_0)\geq \log \frac{\kappa-2}{2c_2(2+\kappa_1)}.$$
	
	$(ii)$ Let $i\in\{1,\cdots,m-1\}$. For any $\lambda$-curve $\beta_i\in \Lambda_{\xi_i\xi_{i+1}}^{\lambda}$, it holds
	$$\ell_k(\beta_i)< 4A_0.$$
	
	$(iii)$ Let $i\in\{1,\cdots,m-1\}$. For any $\lambda$-curve $\beta_i\in \Lambda_{\xi_i\xi_{i+1}}^{\lambda}$ and any $z\in\beta_i$, it holds
	$$k_D(\xi_i,y_0)-4A_0\leq k_D(z,y_0)\leq k_D(\xi_i,y_0)+4A_0.$$
	\elem

	In the following, for each $i\in \{1,\cdots,m\}$, let $\{\xi_i\}$ be the sequence given by \eqref{2022-09-30-3}. Then for each $i\in \{1,\cdots,m-1\}$, we fix a $\lambda$-curve $\beta_i\in \Lambda_{\xi_{i}\xi_{i+1}}^\lambda$ and set $\beta_0=\bigcup_{i=1}^{m-1}\beta_i$; see Figure \ref{fig3}. Fix $u$, $v\in \beta_0$. Then there are $i_1$ and $i_2\in\{1,\cdots,m-1\}$ such that $u\in\beta_{i_1}$ and $v\in\beta_{i_2}$.

	\begin{lem}[{\cite[Lemma 5.12]{Guo-H-W-2025}}]\label{lem-22-09-20}  Suppose that $i_1$, $i_2\in\{1,\cdots,m-1\}$ with $i_1<i_2$, $u\in\beta_{i_1}$ and $v\in\beta_{i_2}$. Fix $\varsigma_1\in \Lambda_{uy_0}^{\lambda}$, $\varsigma_2\in\Lambda_{vy_0}^{\lambda}$ and $\varsigma_3\in\Lambda_{uv}^{\lambda}$.  Then for each $i\in \{1,2,3\}$,  there exists a point $z_{i_1i_2,i}\in \varsigma_{i}$ such that
		\begin{enumerate}
			\item\label{lem-22-20-1}
			For each $j\in\{1,2,3\}$, $$k_D(z_{i_1i_2,j},z_{i_1i_2,j+1})\leq 3C_0,$$ where $z_{i_1i_2,4}=z_{i_1i_2,1}$.
			\item\label{lem-22-20-2}
			$k_D(u,z_{i_1i_2,1})-5A_0< k_D(v,z_{i_1i_2,2})<k_D(u,z_{i_1i_2,1})+5A_0$.
			\item\label{lem-22-09-20-3}
			$\frac{1}{2}k_D(u,v)-\frac{5A_0+6C_0}{2}\leq k_D(u,z_{i_1i_2,3})\leq \frac{1}{2}k_D(u,v)+\frac{1}{8}+\frac{5A_0+6C_0}{2}$.
			\item\label{lem-22-09-20-4}
			$\frac{1}{2}k_D(u,v)-\frac{5A_0+6C_0}{2}\leq k_D(v,z_{i_1i_2,3})\leq \frac{1}{2}k_D(u,v)+\frac{1}{8}+\frac{5A_0+6C_0}{2}$.
		\end{enumerate}
	\end{lem}

	\subsection{Six-tuples with Property \ref{Property-A}}
	In this section, we recall the defintion of six-tuples introduced in \cite[Section 5.5]{Guo-H-W-2025}. Throughout this section, we use $\mathfrak{c}$ to denote a constant with $\mathfrak{c}\geq 64A_0$. 
	Fix $x,y\in D$ and a curve $\beta$ joining $x$ and $y$ in $D$.
	
	\begin{defn}[Six-tuple]\label{def:six-tuple}
		We call $[w_1,w_2,\gamma_{12},\varpi,w_3,w_4]$ a {\it six-tuple}  if the following conditions are satisfied 
		\begin{itemize}
			\item$w_1\in D$ and $w_2\in \beta$;
			\item $\varpi\in \beta[w_2,y]$ and $\gamma_{12}\in \Lambda_{w_1w_2}^\lambda(D)$;
			\item $w_3\in \gamma_{12}$ and $w_4\in \beta[w_2,\varpi]$.
		\end{itemize}
	\end{defn}
	
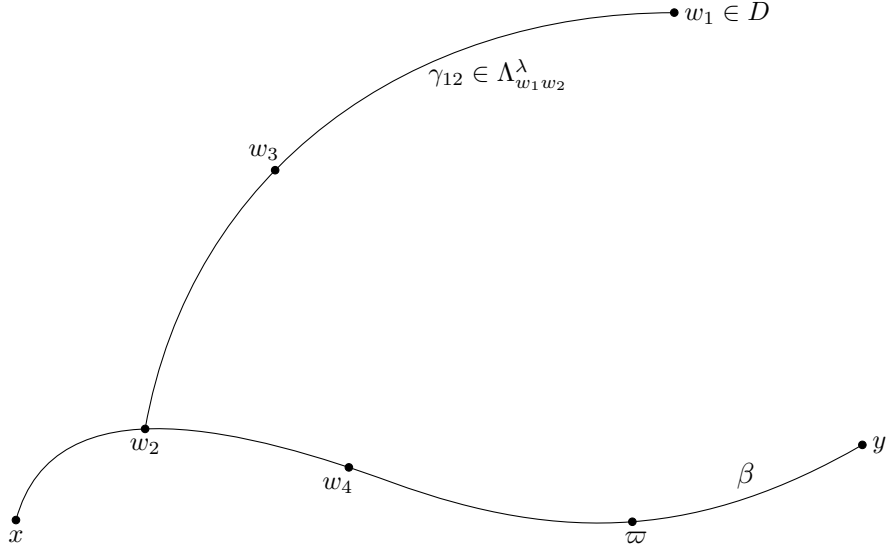
\begin{figure}[htbp]
		\begin{center}
			\begin{tikzpicture}[scale=2.5]
				\draw (-2,-0.2) to [out=75,in=160]coordinate[pos=0] (x) coordinate[pos=0.45] (w2) coordinate[pos=0.9] (w4)  (0,0) to [out=-20,in=210]
				coordinate[pos=0.5] (varpi)
				coordinate[pos=1] (y)
				node [above, pos=0.75]{$\beta$}
				(2.5,0.2);
				\filldraw  (x)node[below] {\small $x$} circle (0.02);
				\filldraw  (w2)node[below] {\small $w_{2}$} circle (0.02);

				\filldraw  (w4)node[below, rotate=0] {\small $w_{4}\;\;\;$} circle (0.02);
				
				\filldraw  (varpi)node[below] {\small $\;\varpi$} circle (0.02);
				\filldraw  (y)node[right] {\small $y$} circle (0.02);
				
				\draw (w2) to [out=80,in=180]
				coordinate[pos=0.4] (w3)
				coordinate[pos=0.75] (gamma12)
				coordinate[pos=1] (w1)
				(1.5,2.5);
				\filldraw  (w3)node[above, xshift=-0.15cm] {\small $w_{3}$} circle (0.02);
				\node[below, xshift=0.15cm, yshift=-0.1cm] at(gamma12) {\small $\gamma_{12}\in \Lambda_{w_1w_2}^{\lambda}$};
				\filldraw  (w1)node[right] {\small $w_{1}\in D$} circle (0.02);

			\end{tikzpicture}
		\end{center}
		\caption{Definition of a six-tuple} \label{fig-8-15-2}
	\end{figure}

	We usually use the symbol $\Omega=[w_1,w_2,\gamma_{12},\varpi,w_3,w_4]$ to represent a six-tuple; see Figure \ref{fig-8-15-2}.
	
	\begin{Prop}\label{Property-A}
		Given $\mathfrak{c}\geq 64A_0$, a six-tuple $\Omega=[w_1,w_2,\gamma_{12},\varpi,w_3,w_4]$ is said to satisfy {\it Property \ref{Property-A}} with constant $\mathfrak{c}$ if
		\begin{enumerate}
			\item\label{3-3.1}
			$k_D(w_1, w_2)\leq \mathfrak{c}$ and $k_D(w_3, w_4)=\mathfrak{c}$;
			\item\label{3-3.2}
			for any $u\in \gamma_{12}$ and any $v\in \beta[w_4,\varpi]$,
			$k_D(u,v)\geq \mathfrak{c}-\frac{1}{4}.$
		\end{enumerate}
		If a six-tuple $\Omega$ satisfies Property \ref{Property-A} with constant $\mathfrak{c}$, then we shall write $\Omega^{\mathfrak{c}}=\Omega$ to indicate its dependence on the constant $\mathfrak{c}$.
	\end{Prop}
	In most situations, we shall simply say that $\Omega^{\mathfrak{c}}=[w_1,w_2,\gamma_{12},\varpi,w_3,w_4]$ is a six-tuple with Property \ref{Property-A} to indicate that the six-tuple $[w_1,w_2,\gamma_{12},\varpi,w_3,w_4]$ satisfies Property \ref{Property-A} with constant $\mathfrak{c}$.

	The following result provides a sufficient condition for the existence of a six-tuples with Property \ref{Property-A}. 
	\begin{lem}[{\cite[Lemma 5.14]{Guo-H-W-2025}}]\label{lem-3.8}
		Let $w_1\in D$, $w_2\in\beta$ and $\gamma_{12}\in \Lambda_{w_1w_2}^\lambda(D)$. Suppose that $0<k_D(w_1,w_2)\leq \mathfrak{c}$ and there is $\varpi\in\beta[w_2,y]$ such that for any $w\in \gamma_{12}$,
		$k_D(w,\varpi)>\mathfrak{c}$.
		Then there exist two points $w_3\in \gamma_{12}\backslash\{w_2\}$ and $w_4\in\beta[w_2,\varpi]$
		such that $\Omega^{\mathfrak{c}}=[w_1,w_2,\gamma_{12},\varpi,w_3,w_4]$ is a six-tuple with Property \ref{Property-A}.
	\end{lem}

	The next result gives elementary estimates for quasihyperbolic distance of points related to a six-tuple with Property \ref{Property-A}.
	\begin{lem}[{\cite[Lemma 5.16]{Guo-H-W-2025}}]\label{lem23-01} 
		Suppose that $[w_1,w_2,\gamma_{12},\varpi,w_3,w_4]$ is a six-tuple satisfying {Property \ref{Property-A}} with constant $\mathfrak{c}\geq 36C_0$. For each integer $r\in \{1,2,3\}$ and $s\in\{r+1,\cdots,4\}$, let $\gamma_{rs}\in \Lambda_{w_rw_s}^\lambda(D)$, where $\gamma_{13}=\gamma_{12}[w_1,w_3]$ and $\gamma_{23}=\gamma_{12}[w_2,w_3]$. Then
		\begin{enumerate}
			\item\label{2022-08-20}
			\begin{enumerate}
				\item\label{23-07-26-4}
				there exist $x_1\in \gamma_{23}$, $x_2\in\gamma_{24}$ and $x_3\in\gamma_{34}$ such that for any $r\not=s\in\{1,2,3\}$, $$k_D(x_r,x_s)\leq 3C_0;$$ 
				\item\label{23-07-26-5}
				$k_D(w_3,x_3)\leq \frac{1}{2}+3C_0$.
			\end{enumerate}
			
			\item\label{2022-08-30}
			\begin{enumerate} \item\label{23-07-26-6} there exist $p_1\in \gamma_{13}$, $p_2\in\gamma_{14}$ and $p_3\in\gamma_{34}$ such that for any $r\not=s\in\{1,2,3\}$, $$k_D(p_r,p_s)\leq 3C_0;$$ 
				\item\label{23-07-26-6} $k_D(w_3,p_3)\leq \frac{1}{2}+3C_0$.
			\end{enumerate}
			
			\item\label{23-01-10-0} $k_D(w_2,w_4)\geq  k_D(w_2,w_3)+k_D(w_3,w_4)-\frac{5}{4}-15C_0$.
			
			\item\label{23-12-1-Add} for any $w\in\gamma_{34}$, it holds
			$$k_D(y_0,w)\leq k_D(y_0,\xi_1)-\min\{k_D(w_2,w_3)+k_D(w_3,w), k_D(w_4,w)\}+14A_0.$$
		\end{enumerate}
	\end{lem}

	The following lemma gives estimate on the quasihyperbolic distance between $w_1$ and $w_3$.
	\begin{lem}[{\cite[Lemma 5.15 and Lemma 5.19]{Guo-H-W-2025}}]\label{lem-23-01} 
		Let $\gamma_1$ be given as in Section \ref{subsec:xi sequence}  $($see Figure \ref{fig3}$)$. Let $w_{1,1}\in D$, $w_{1,2}\in \beta_0$ and $\gamma_{12}=\gamma_1[w_{1,1},w_{1,2}]$, $w_{1,3}\in \gamma_{12}$, $\varpi\in\beta_0[w_{1,2},\xi_m]$ and $w_{1,4}\in \beta_0[w_{1,2},\varpi]$. 
		\ben
		\item\label{H25-086-1}  If $w_{1,1}\in\gamma_1$, $w_{1,2}=\xi_1$ and the six-tuple $\Omega_1=[w_{1,1},w_{1,2},\gamma_{12},\varpi,w_{1,3},w_{1,4}]$ satisfies Property \ref{Property-A} with constant $\mathfrak{c}\geq 16A_0$. Then 
		$$k_D(w_{1,1},w_{1,3})\leq 11A_0.$$	
		
		\item\label{H25-086-1} Let $w_{2,1}\in\gamma_{12}[w_{1,2},w_{1,3}]$, $w_{2,2}=w_{1,2}$, $\gamma_{12}^1=\gamma_{12}[w_{2,2},w_{2,1}]$, $\varpi_2=w_{1,4}\in \beta_0[w_{2,2},\xi_m]$, $w_{2,3}\in \gamma_{12}^1$ and $w_{2,4}\in \beta_0[w_{2,2},\varpi_2]$. If the six-tuple $\Omega_2^{\mathfrak{c}}=[w_{2,1},w_{2,2},\gamma_{12}^1,\varpi_2,w_{2,3},w_{2,4}]$ also satisfies Property \Ref{Property-A} and $k_D(w_{1,1},w_{1,2})=\mathfrak{c}$, then 
		$$k_D(w_{2,1},w_{2,3})\leq 11A_0.$$
		\een    
	\end{lem}
	\medskip

	Let $w_{1,1}\in D$, $\varpi\in \beta_0$, $w_{1,2}\in \beta_0[\xi_1, \varpi]$, $\gamma_{12}^1\in \Lambda_{w_{1,1}w_{1,2}}^{\lambda}(D)$ $w_{1,3}\in \gamma_{12}^1$ and $w_{1,4}\in \beta_0[w_{1,2},\varpi]$. Then $[w_{1,1},w_{1,2},\gamma_{12}^1,\varpi,w_{1,3},w_{1,4}]$ is a six-tuple. We let $w_{2,1}=w_{1,3}$, $w_{2,2}=w_{1,4}$, $\gamma_{12}^2\in \Lambda_{w_{2,1}w_{2,2}}^{\lambda}(D)$, 
	$w_{2,3}\in \gamma_{12}^{2}$ and $w_{2,4}\in \beta_0[w_{2,2}, \varpi]$. Then $[w_{2,1},w_{2,2},\gamma_{12}^2,\varpi,w_{2,3},w_{2,4}]$ is a six-tuple as well. 
	For $i\in \{1,2\}$, we set $$\Omega_i=[w_{i,1},w_{i,2},\gamma_{12}^i,\varpi, w_{i,3},w_{i,4}].$$
	
	The next lemma gives an upper estimate on $k_D(w_{1,3},w_{2,3})$.
	\begin{lem}[{\cite[Lemma 5.17]{Guo-H-W-2025}}]\label{lem-22-12-0} 
		If both $\Omega_1^{\mathfrak{c}}$ and $\Omega_2^{\mathfrak{c}}$ satisfy Property \ref{Property-A} with constant $\mathfrak{c}\geq 36A_0$, then
		$$k_D(w_{1,3},w_{2,3})=k_D(w_{2,1},w_{2,3})\leq \frac{1}{2}(\mathfrak{c}-k_D(w_{1,2},w_{1,3}))+28A_0.$$
	\end{lem}
	\medskip

	\subsection{A sequence of new six-tuples}
	
	In the rest of this paper, we make the following notational convention:
	$$\rho_1>0\;\;\;\mbox{and }\;\;\frac{\rho_3}{4}\geq\rho_2>\max\{\rho_1,e^{36A_0}\}.$$

	
	
	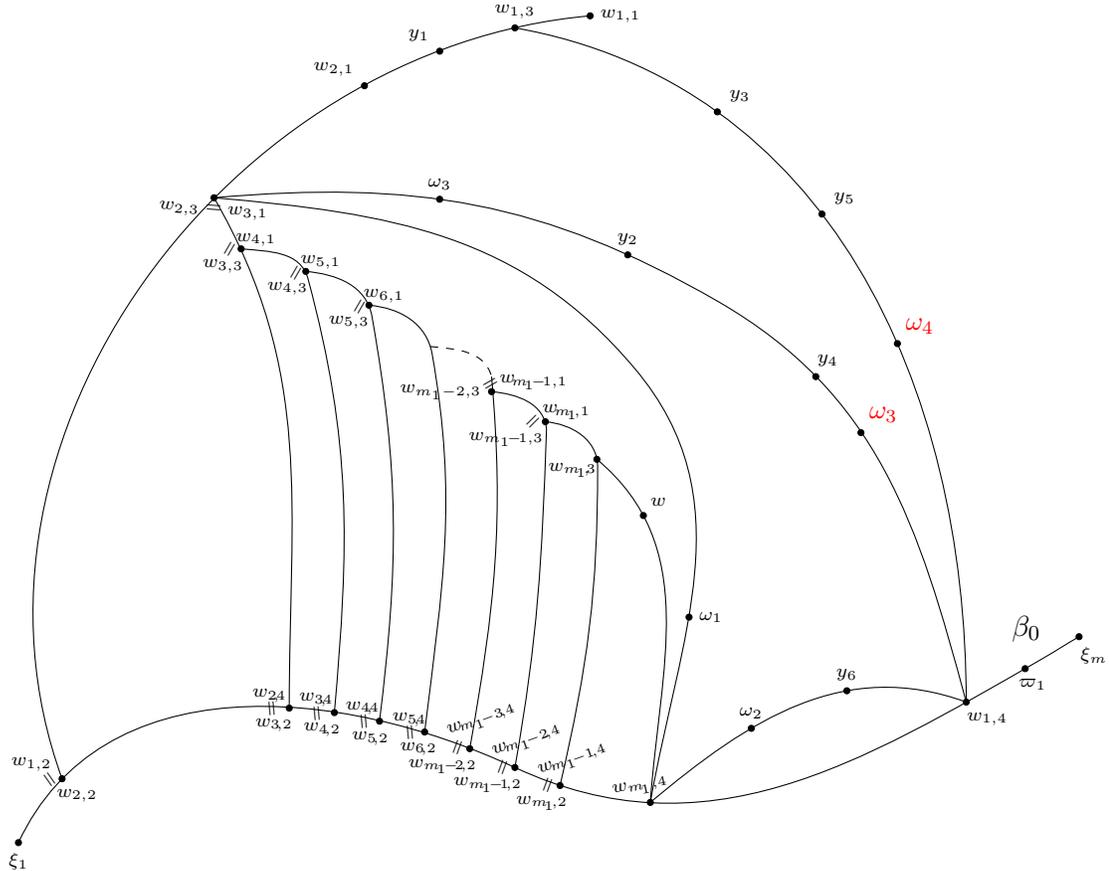
\begin{figure}[htbp]
		\begin{center}
			\begin{tikzpicture}[scale=1]
				\draw (-6.6,-1) to [out=65,in=155] coordinate[pos=0] (xi0) (0,0) to [out=-25,in=210] coordinate[pos=0.97] (varpi) (7,1.44) to [out=30,in=210]  node[above,pos=0,xshift=-0.2cm,yshift=0.1cm] {$\beta_0$} coordinate[pos=1] (xim) (7.5,1.74);
				\filldraw  (xi0)node[below] {\tiny $\xi_1$} circle (0.04);
				\filldraw  (-6.02,-0.15)node[below] {\tiny $\;\;\;\;\;w_{2,2}$} circle (0.04);
				\node[left] at(-6.02,0.05) {\tiny $w_{1,2}$};
				\node[rotate=-55] at(-6.2,-0.17) {\tiny $=$};

				\filldraw  (-3,0.79)node[below] {\tiny $w_{3\!,2}\;\;\;\;\;$} circle (0.04);
				\node[above] at(-3,0.74) {\tiny $w_{2\!,\!4}\;\;\;\;\;\;\,$};
				\node[rotate=-85] at(-3.25,0.78) {\tiny $=$};
				
				\filldraw  (-2.4,0.73)node[below] {\tiny $w_{4\!,2}\;\;\;\;$} circle (0.04);
				\node[above] at(-2.4,0.68) {\tiny $w_{3\!,\!4}\;\;\;\;\;\;$};
				\node[rotate=-85] at(-2.65,0.72) {\tiny $=$};
				
				\filldraw  (-1.8,0.616)node[below] {\tiny $w_{5\!,2}\;\;\;$} circle (0.04);
				\node[above] at(-1.8,0.57) {\tiny $w_{4\!,\!4}\;\;\;\;\;\,$};
				\node[rotate=-85] at(-2.02,0.616) {\tiny $=$};
				
				\filldraw  (-1.2,0.472)node[below] {\tiny $w_{6\!,2}\;\;$} circle (0.04);
				\node[above] at(-1.2,0.422) {\tiny $w_{5\!,\!4}\;\;\;\;\;$};
				\node[rotate=-85] at(-1.41,0.472) {\tiny $=$};
				
				\filldraw  (-0.6,0.25)node[below] {\tiny $w_{m_{\!1}\!-\!2\!,2}\;\;\;\;\;\;\;\;\;$} circle (0.04);
				\node[rotate=70] at(-0.75,0.27) {\tiny $=$};
				\node[rotate=20] at(-0.44,0.63) {\tiny $w_{\!m_{1}\!-3\!,4}$};
				
				\filldraw  (0,0)node[below] {\tiny $w_{m_{\!1}\!-\!1\!,2}\;\;\;\;\;\;\;\;\;$} circle (0.04);
				\node[rotate=70] at(-0.15,0) {\tiny $=$};
				\node[rotate=20] at(0.16,0.35) {\tiny $w_{\!m_{1}\!-2\!,4}$};
				
				\filldraw  (0.6,-0.24)node[below] {\tiny $w_{m_{\!1}\!,2}\;\;\;\;\;\;$} circle (0.04);
				\node[rotate=70] at(0.45,-0.23) {\tiny $=$};
				\node[rotate=20] at(0.77,0.1) {\tiny $w_{\!m_{1}\!-1\!,4}$};
				
				\filldraw  (1.8,-0.465) circle (0.04);
				\node[rotate=12] at(1.59,-0.27) {\tiny $\;\;w_{m_{\!1}\:\!,4}$};

				\filldraw  (6,0.87)node[below] {\tiny $\;\;\;\;\;\;\;w_{1,4}$} circle (0.04);
				
				\filldraw  (varpi)node[below,xshift=0.1cm] {\tiny $\varpi_1$} circle (0.04);
				\filldraw  (xim)node[below,xshift=0.2cm] {\tiny $\xi_m$} circle (0.04);

				\filldraw  (1,10)node[right] {\tiny $w_{1,1}$} circle (0.04);
				\draw (-6.02,-0.15)  to [out=110,in=185] (1,10);
				\filldraw  (0,9.84)node[above] {\tiny $w_{1,3}$} circle (0.04);
				\filldraw  (-1,9.53)node[above left] {\tiny $y_1$} circle (0.04);
				\filldraw  (-2,9.07)node[above left] {\tiny $w_{2,1}$} circle (0.04);
				
				\filldraw  (-4,7.58) circle (0.04);
				\node[below, rotate=-4] at(-4.3,7.65) {\tiny $w_{2,3}=$};
				\node[below] at(-3.55,7.6) {\tiny $w_{3,1}$};
				
				\draw (-4,7.58)  to [out=-60,in=90] (-3,0.79);
				\filldraw  (-3.64,6.9)node[below] {\tiny $w_{3,3}\;\;\;\;\;\;$} circle (0.04);
				\node[rotate=60] at(-3.78,6.89) {\tiny $=$};
				\node[above] at(-3.64,6.8) {\tiny $\;\;\;\;\;w_{4,1}$};
				
				\draw (-3.64,6.9) to [out=-5,in=105] (-2.75,6.5)  to [out=-75,in=85] (-2.4,0.73);
				\filldraw  (-2.78,6.6)node[below] {\tiny $w_{4,3}\;\;\;\;\;\;$} circle (0.04);
				\node[rotate=60] at(-2.9,6.58) {\tiny $=$};
				\node[above] at(-2.78,6.5) {\tiny $\;\;\;\;\;w_{5,1}$};
				
				\draw (-2.78,6.6) to [out=-5,in=100] (-1.9,6)  to [out=-80,in=83] (-1.8,0.616);
				\filldraw  (-1.94,6.15)node[below] {\tiny $w_{5,3}\;\;\;\;\;\;\,$} circle (0.04);
				\node[rotate=60] at(-2.05,6.13) {\tiny $=$};
				\node[above] at(-1.94,6.05) {\tiny $\;\;\;\;\;w_{6,1}$};
				
				\draw (-1.94,6.15) to [out=-5,in=100] (-1.1,5.5)  to [out=-80,in=83] (-1.2,0.472);
				
				\draw[dashed] (-1.125,5.6) to [out=-5,in=105] (-0.31,5.2);
				\draw (-0.31,5.2)  to [out=-85,in=80] (-0.6,0.25);
				\filldraw  (-0.31,5)node[left] {\tiny $w_{m_{1}-2,3}$} circle (0.04);
				\node[rotate=30] at(-0.32,5.1) {\tiny $=$};
				\node[above] at(-0.31,4.9) {\tiny $\;\;\;\;\;\;\;\;\;\;\;\;\;\;w_{\!m_{\!1}\!-\!1,1}$};
				
				\draw (-0.31,5) to [out=-5,in=95] (0.42,4.5)  to [out=-92,in=80] (0,0);
				\filldraw  (0.405,4.6)node[below] {\tiny $w_{m_{1}\!-\!1,3}\;\;\;\;\;\;\;\;\;\;\;\;\;$} circle (0.04);
				\node[rotate=45] at(0.25,4.6) {\tiny $=$};
				\node[above] at(0.405,4.5) {\tiny $\;\;\;\;\;\;\;w_{\!m_{\!1}\!,1}$};
				
				\draw (0.405,4.6) to [out=-5,in=90] (1.1,4)  to [out=-90,in=77] (0.6,-0.24);
				\filldraw  (1.09,4.1) circle (0.04);
				\node[below] at(0.767,4.19) {\tiny $w_{m_{\!1}\!,\!3}$};
				
				\draw (1.09,4.1) to [out=-40,in=86] coordinate[pos=0.2] (w) (1.8,-0.465);
				\filldraw  (w)node[above,xshift=0.2cm] {\tiny$w$} circle (0.04);
				
				\draw (-4,7.58) to [out=-5,in=130] (1.6,5.3) to [out=-50,in=80] (1.8,-0.465);
				\filldraw  (2.315,2)node[right] {\tiny $\omega_1$} circle (0.04);
				
				\draw (-4,7.58) to [out=5,in=150] (2.7,6.2) to [out=-30,in=105] (6,0.87);
				\filldraw  (-1,7.56)node[above] {\tiny $\omega_3$} circle (0.04);
				\filldraw  (1.5,6.82)node[above] {\tiny $y_2$} circle (0.04);
				\filldraw  (4,5.2)node[above] {\tiny $\;\;\;\;y_4$} circle (0.04);
				\filldraw  (4.6,4.455)node[above,xshift=0.3cm] {\small ${\color{red} \omega_3}$} circle (0.04);
				
				\draw (0,9.84) to [out=-10,in=90] coordinate[pos=0.6] (w4) coordinate[pos=0.425] (eta2) coordinate[pos=0.25] (y3) (6,0.87);
				\filldraw  (w4)node[above,xshift=0.3cm] {\small ${\color{red} \omega_4}$} circle (0.04);
				\filldraw  (eta2)node[above,xshift=0.3cm] {\tiny$y_5$} circle (0.04);
				\filldraw  (y3)node[above,xshift=0.3cm] {\tiny $y_3$} circle (0.04);
				
				\draw (1.8,-0.465) to [out=40,in=160] coordinate[pos=0.34] (w2) coordinate[pos=0.65] (eta3)  (6,0.87);
				\filldraw  (w2)node[above] {\tiny $\omega_2$} circle (0.04);
				\filldraw  (eta3)node[above] {\tiny $y_6$} circle (0.04);

			\end{tikzpicture}
		\end{center}
		\caption{Illustration for the construction of new six-tuples} \label{fig-8-16-2}
	\end{figure}

	We are going to construct a finite sequence of six-tuples. Let $\Omega_1=[w_{1,1},w_{1,2},\gamma_{12}^1,\varpi_1,w_{1,3},w_{1,4}]$ be a six-tuple, where $w_{1,1}\in D$, $w_{1,2}\in \beta_0$, $\gamma_{12}^1\in \Lambda_{w_{1,1}w_{1,2}}^{\lambda}$, 
	$\varpi_1\in\beta_0[w_{1,2},\xi_m]$, $w_{1,3}\in \gamma_{12}^1$ and $w_{1,4}\in \beta_0[w_{1,2},\varpi_1]$ (see Figure \ref{fig-8-16-2}). Suppose
	\be\label{eq:13 rho 1}
	k_D(w_{1,1},w_{1,2})=\rho_3\quad\text{and}\quad k_D(w_{1,1},w_{1,3})=\rho_1.
	\ee
	Then 
	\[
	k_{D}(w_{1,2},w_{1,3})\geq k_{D}(w_{1,2},w_{1,1})-k_{D}(w_{1,1},w_{1,3})\geq \rho_3-\rho_1\geq 2\rho_2.
	\]
	This implies that there exists some $w_{2,1}\in \gamma_{12}^1[w_{1,2},w_{1,3}]$ such that
	\be\label{23-07-29-10}
	k_D(w_{1,3},w_{2,1})=\rho_2,
	\ee
	
	Based on $\Omega_1$ and $w_{2,1}$, we obtain the second six-tuple $\Omega_2=[w_{2,1},w_{2,2},\gamma_{12}^2,\varpi_2,w_{2,3},w_{2,4}]$ by setting $w_{2,2}=w_{1,2}$, $\gamma_{12}^2=\gamma_{12}^1[w_{2,1},w_{2,2}]$, $\varpi_2=w_{1,4}\in \beta_0[w_{2,2},\xi_m]$, $w_{2,3}\in \gamma_{12}^2$ and $w_{2,4}\in \beta_0[w_{2,2},\varpi_2]$ (see Figure \ref{fig-8-16-2}).
	
	Then we obtain the third six-tuple $\Omega_3=[w_{3,1},w_{3,2},\gamma_{12}^3,\varpi_3,w_{3,3},w_{3,4}]$ by setting $w_{3,1}=w_{2,3}$, $w_{3,2}=w_{2,4}$, $\gamma_{12}^3\in \Lambda_{w_{2,3}w_{2,4}}^\lambda$, $\varpi_3=w_{1,4}\in \beta_0[w_{3,2},\xi_m]$, $w_{3,3}\in \gamma_{12}^3$ and $w_{3,4}\in \beta_0[w_{3,2},\varpi_3]$
	and the fourth six-tuple $\Omega_4=[w_{4,1},w_{4,2},\gamma_{12}^4,\varpi_4,w_{4,3},w_{4,4}]$ by setting $w_{4,1}=w_{3,3}$, $w_{4,2}=w_{3,4}$, $\gamma_{12}^4\in \Lambda_{w_{3,3}w_{3,4}}^\lambda$, $\varpi_4=w_{1,4}\in \beta_0[w_{4,2},\xi_m]$, 
	$w_{4,3}\in \gamma_{12}^{4}$ and $w_{4,4}\in \beta_0[w_{4,2},\varpi_4]$ (see Figure \ref{fig-8-16-2}).
	
	By repeating this procedure, we obtain a finite sequence of six-tuples $\{\Omega_i\}_{i=1}^{m_1}$ with $\Omega_i=[w_{i,1},w_{i,2},\gamma_{12}^i,\varpi_i,w_{i,3},w_{i,4}]$ by setting 
	for each $i\in\{3,\cdots,m_1\}$, $w_{i,1}=w_{i-1,3}$, $w_{i,2}=w_{i-1,4}$, $\gamma_{12}^i\in \Lambda_{w_{i-1,3}w_{i-1,4}}^\lambda$, $\varpi_i=w_{1,4}\in \beta_0[w_{i,2},\xi_m]$, $w_{i,3}\in \gamma_{12}^i$ and $w_{i,4}\in \beta_0[w_{i,2},\varpi_i]$.

	%
	%
	%
	
	Then we have the following useful estimates.
	
	\begin{lem}\label{lem-22-09-30} 
		$(1)$
		Suppose that $m_1\geq 2$ and $\Omega_1^{\rho_3}$ satisfies Property \ref{Property-A}. Then	
		\beqq
		k_D(w_{1,4},w_{2,3})\geq k_D(w_{1,3},w_{1,4})+k_D(w_{1,3},w_{2,3})-\frac{5}{4}-12C_0.
		\eeqq
		
		\noindent $(2)$
		Suppose that for each $i\in\{1,\cdots,m_1\}$, $\Omega_i^{\rho_3}$ satisfies Property \Ref{Property-A}.
		Then 
		\beqq
		k_D(w_{i,3},w_{i+1,3})\leq \frac{\rho_1+\rho_2}{2^{i-1}}+68A_0 \qquad\quad \text{for all }i\in\{2,\cdots,m_1-1\}.
		\eeqq
		
		\noindent $(3)$
		Suppose $m_1\geq 3$ and for each $i\in\{1,\cdots,m_1\}$, $\Omega_i^{\rho_3}$ satisfies Property \ref{Property-A}. If $k_D(w_{2,3},w_{m_1,3})\leq \rho_2-\rho_1-24A_0e^{36A_0}$, then 
		
		$(3a)$ $k_D(w_{1,4},w_{2,3})> k_D(w_{2,3},w_{m_1,4})-\frac{7}{4}+\rho_1+24A_0e^{36A_0}$, and
		
		$(3b)$
		for any $\gamma_{34}^{m_1}\in\Lambda_{w_{m_1,3}w_{m_1,4}}^{\lambda}$ and $w\in \gamma_{34}^{m_1}$, it holds 
		$
		k_D(w_{1,4},w)> \rho_3.$
		
	\end{lem}
	
	
	\begin{proof}
		Since statements (1) and (3) follow from \cite[Lemma 5.18(1)]{Guo-H-W-2025} and  \cite[Lemma 5.18(3)]{Guo-H-W-2025}, we only need to prove (2).
		
		It follows from Lemma \ref{lem-22-12-0} and the above construction that 
		\be\label{H25-0810-1}k_D(w_{2,3},w_{3,3})\leq \frac{1}{2}(\rho_3-k_D(w_{2,2},w_{2,3}))+28A_0\leq \frac{1}{2}k_D(w_{1,1},w_{2,3})+28A_0.\ee
		
		For each $i\in\{2,\cdots,m_1-1\}$, \eqref{eq:13 rho 1} and our assumption in the lemma imply $$k_D(w_{i-1,3},w_{i,3})=k_D(w_{i,1},w_{i,3})\geq k_D(w_{i,1},w_{i,2})-k_D(w_{i,2},w_{i,3})\geq \rho_3-k_D(w_{i,2},w_{i,3}).$$
		Thus, for each $i\in\{3,\cdots,m_1-1\}$, we may repeatedly apply Lemma \ref{lem-22-12-0} to obtain 
		\begin{equation}\label{eq:inductive estimate}
			\begin{aligned}
				k_D(w_{i,3},w_{i+1,3})&\leq \frac{1}{2}(\rho_3-k_D(w_{i,2},w_{i,3})+28A_0\leq \frac{1}{2}k_D(w_{i-1,3},w_{i,3})+28A_0\\
				&\leq \frac{1}{2^{i-2}}k_D(w_{2,3},w_{3,3})+28A_0\left(1+\frac{1}{2}+\cdots+\frac{1}{2^{i-3}}\right)\\&\stackrel{\eqref{H25-0810-1}}{\leq}  \frac{1}{2^{i-1}}k_D(w_{1,1},w_{2,3})+28A_0\left(1+\frac{1}{2}+\cdots+\frac{1}{2^{i-2}}\right)\\
				&< \frac{1}{2^{i-1}}k_D(w_{1,1},w_{2,3})+ 56A_0.	
			\end{aligned}
		\end{equation}
		Then it follows from \eqref{H25-0810-1} and \eqref{eq:inductive estimate} that for each $i\in\{2,\cdots,m_1-1\}$,
		\[
		\begin{aligned}
			k_D(w_{i,3},w_{i+1,3})&< \frac{1}{2^{i-1}}k_D(w_{1,1},w_{2,3})+56A_0\\
			&\leq \frac{k_D(w_{1,1},w_{1,3})+k_D(w_{1,3},w_{2,1})+k_D(w_{2,1},w_{2,3})}{2^{i-1}}+56A_0\\
			&\stackrel{\eqref{eq:13 rho 1} + \eqref{23-07-29-10}}{\leq} \frac{\rho_1+\rho_2}{2^{i-1}}+\frac{1}{2^{i-1}}k_D(w_{2,1},w_{2,3})+56A_0\stackrel{\text{Lemma}\; \ref{lem-23-01}}{\leq} \frac{\rho_1+\rho_2}{2^{i-1}}+68A_0.
		\end{aligned}
		\]
		This proves (2). 
		
		
		
		
		
	\end{proof}
	%
	%
	%

	Let $w_{0,1}\in D$, $w_{0,2}\in \beta_0[\xi_1,\xi_m]$, $\gamma_{12}^0\in \Lambda_{w_{0,1},w_{0,2}}^\lambda$,  $\varpi_0\in\beta_0[\xi_{1},\xi_m]$, $w_{0,3}\in \gamma_{12}^0$ and $w_{0,4}\in \beta[w_{0,2},\varpi_0]$. Suppose $\Omega_0^{\rho_{3}}=[w_{0,1},w_{0,2},\gamma_{12}^{0},\varpi_0,w_{0,3},w_{0,4}]$ satisfies Property \Ref{Property-A},
	\be\label{24-2-15-1}
	k_D(w_{0,1},w_{0,2})=\rho_{3}.
	\ee
	and
	\be\label{24-2-26-10}
	k_D(w_{0,1},w_{0,3})=\rho_{1}>0
	\ee
	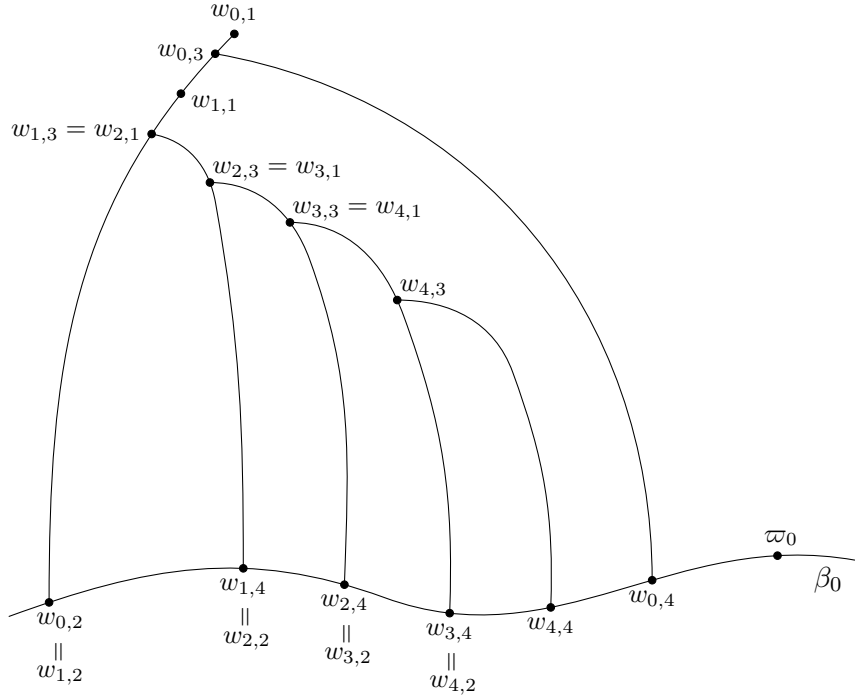
\begin{figure}[htbp]
		\begin{center}
			\begin{tikzpicture}[scale=2.5]
				\draw (-2,-0.1) to [out=20,in=160]
				coordinate[pos=0.1] (w02)
				coordinate[pos=0.63] (w14)
				coordinate[pos=0.9] (w24)
				(0,0) to [out=-20,in=170]
				coordinate[pos=0.13] (w34)
				coordinate[pos=0.35] (w44)
				coordinate[pos=0.58] (w04)
				coordinate[pos=0.85] (varpi)
				node [below, pos=0.95]{$\beta_0$}
				(2.5,0.2);
				\filldraw  (w02)node[below, xshift=0.15cm] {\small $w_{0,2}$} circle (0.02);
				\node[below,rotate=90, xshift=-0.65cm, yshift=0.05cm] at(w02) {\small $=$};
				\node[below, xshift=0.15cm, yshift=-0.7cm] at(w02) {\small $w_{1,2}$};

				\filldraw  (w14)node[below, rotate=0] {\small $w_{1,4}$} circle (0.02);
				\node[below,rotate=90, xshift=-0.65cm, yshift=0.15cm] at(w14) {\small $=$};
				\node[below, xshift=0.05cm, yshift=-0.7cm] at(w14) {\small $w_{2,2}$};
				
				\filldraw  (w24)node[below] {\small $w_{2,4}$} circle (0.02);
				\node[below,rotate=90, xshift=-0.65cm, yshift=0.15cm] at(w24) {\small $=$};
				\node[below, xshift=0.05cm, yshift=-0.7cm] at(w24) {\small $w_{3,2}$};
				
				\filldraw  (w34)node[below] {\small $w_{3,4}$} circle (0.02);
				\node[below,rotate=90, xshift=-0.65cm, yshift=0.15cm] at(w34) {\small $=$};
				\node[below, xshift=0.05cm, yshift=-0.7cm] at(w34) {\small $w_{4,2}$};
				
				\filldraw  (w44)node[below] {\small $w_{4,4}$} circle (0.02);
				\filldraw  (w04)node[below] {\small $w_{0,4}$} circle (0.02);
				\filldraw  (varpi)node[above] {\small $\;\varpi_0$} circle (0.02);
				
				\draw (w02) to [out=90,in=225]
				coordinate[pos=0.8] (w13)
				coordinate[pos=0.88] (w11)
				coordinate[pos=0.96] (w03)
				coordinate[pos=1] (w01)
				(-0.8,3);
				
				\filldraw  (w13)node[left,rotate=0] {\small $w_{1,3}=w_{2,1}$} circle (0.02);
				\filldraw  (w03)node[left] {\small $w_{0,3}$} circle (0.02);
				\filldraw  (w11)node[right, yshift=-0.15cm] {\small $w_{1,1}$} circle (0.02);
				\filldraw  (w01)node[above] {\small $w_{0,1}$} circle (0.02);
				
				\draw (w14) to [out=90,in=-80] (-0.9,2.1) to [out=100,in=-10] coordinate[pos=0.2] (w23) (w13);
				
				\filldraw  (w23)node[right,xshift=-0.1cm, yshift=0.15cm, rotate=0] {\small $w_{2,3}=w_{3,1}$} circle (0.02);
				
				\draw (w24) to [out=88,in=-70] (-0.4,1.8) to [out=110,in=0] coordinate[pos=0.3] (w33) (w23);
				
				\filldraw  (w33)node[right,xshift=-0.1cm, yshift=0.15cm, rotate=0] {\small $w_{3,3}=w_{4,1}$} circle (0.02);
				
				\draw (w34) to [out=88,in=-70] (0.1,1.5) to [out=110,in=0] coordinate[pos=0.1] (w43) (w33);
				
				\filldraw  (w43)node[right,xshift=-0.1cm, yshift=0.15cm, rotate=0] {\small $w_{4,3}$} circle (0.02);
				
				\draw (w44) to [out=88,in=-70] (0.68,1.2) to [out=110,in=0] (w43);
				
				\draw (w04) to [out=90,in=-10] (w03);

			\end{tikzpicture}
		\end{center}
		\caption{Illustration for the proof of Lemma \ref{lem-24-05}} \label{fig-8-15-1}
	\end{figure}

	Then we have the following key result (see Figure \ref{fig-8-15-1}). 
	\begin{lem}\label{lem-24-05}
		Suppose that $\rho_{2}\geq 5\rho_1+24A_0e^{36A_0}$. Then
		we can find four six-tuples $\{\Omega_i\}_{i=1}^{4}$ with $\Omega_i^{\rho_3}=[w_{i,1},w_{i,2},\gamma_{12}^i,\varpi_i,w_{i,3},w_{i,4}]$ such that
		for each $i\in\{2,\cdots,4\}$, $\Omega_i^{\rho_3}$ satisfies Property \Ref{Property-A}, where $w_{i,1}=w_{i-1,3}$ and $w_{i,2}=w_{i-1,4}$. Furthermore, we have
		$$k_D(w_{3,3},w_{4,3})\leq \frac{\rho_{1}+\rho_{2}}{2^3}+68A_0\quad\mbox{ and }\quad k_D(w_{0,2},w_{4,2})\geq \rho_3-\frac{1}{4}.$$
	\end{lem}
	
	\bpf 
	Let $w_{1,2}=w_{0,2}$ and $\varpi_1=w_{0,4}$. Then we have the following claim: 
	\medskip
	
	\noindent {\bf Claim \ref{lem-24-05}-{\rm I}}: There are $w_{1,1}\in\gamma_{12}^{0}[w_{0,3},w_{1,2}]$, $\gamma_{12}^1\in \Lambda_{w_{1,1}w_{1,2}}^\lambda$, $w_{1,3}\in \gamma_{12}^{1}$ and $w_{1,4}\in \beta_0[w_{1,2},\varpi_1]$ $($see Figure \ref{fig-8-15-1}$)$ such that
\begin{enumerate}
	\item\label{23-08-29-3} {\it
		the six-tuple $\Omega_1^{\rho_3}=[w_{1,1},w_{1,2},\gamma_{12}^1,\varpi_1,w_{1,3},w_{1,4}]$ satisfies Property \ref{Property-A};  }
	\item\label{23-08-29-4} {\it
		$k_D(w_{0,3},w_{1,1})=\rho_2$, and   }
	\item\label{23-08-29-5} {\it
		for any $u\in \gamma_{12}^{1}$ and any $v\in \beta_0[w_{1,4},\varpi_1]$,
		$k_D(u,v)\geq \rho_3-\frac{1}{4}$.}
\end{enumerate}
\smallskip

\noindent {\bf Proof of Claim \ref{lem-24-05}-{\rm I}}. We first find the point $w_{1,1}$. Since 
\be\label{8-27-1}
k_D(w_{0,3},w_{1,2})\geq k_D(w_{1,2},w_{0,1})-k_D(w_{0,1},w_{0,3}) \stackrel{\eqref{24-2-15-1}+\eqref{24-2-26-10}}{=} \rho_3-\rho_{1}>3\rho_2,
\ee
there exists $w_{1,1}\in\gamma_{12}^{0}[w_{0,3},w_{1,2}]$ $($see Figure \ref{fig-8-15-1}$)$ such that
\be\label{23-08-09-1}k_D(w_{0,3},w_{1,1})=\rho_2.\ee
Then, we know from Lemma \ref{lem-2.2} that
\beq\label{8-26-1}
\begin{aligned}
	k_D(w_{1,1},w_{1,2})  &\leq  \ell_k(\gamma_{12}^{0}[w_{1,1},w_{1,2}] < \ell_k(\gamma_{12}^{0})- \ell_k(\gamma_{12}^{0}[w_{0,3},w_{1,1}])\\  &\leq
	k_D(w_{0,1},w_{0,2})-k_D(w_{0,3},w_{1,1})-k_D(w_{0,1},w_{0,3})+\frac{1}{4}\\  &\stackrel{\eqref{24-2-15-1}+\eqref{23-08-09-1}}{=}
	\frac{1}{4}-\rho_1-\rho_2 +\rho_3.
\end{aligned}
\eeq

Set $\gamma_{12}^1=\gamma_{12}^0[w_{1,1},w_{1,2}]$ and we shall apply Lemma \ref{lem-3.8} to find the remaining two points $w_{1,3}$ and $w_{1,4}$. For this, we need to show that for any $v\in \gamma_{12}^{1}$, $\rho_3$ is a lower bound for the quantity $k_D(v,w_{0,4})$. 

Since the six-tuple $\Omega_0^{\rho_3}$ satisfies Property \ref{Property-A}, we see from
Lemma \ref{lem23-01}(\ref{23-01-10-0}) that
\be\label{23-08-20-2}
\begin{aligned}
	k_D(w_{0,4},w_{1,2})&=  k_D(w_{0,2},w_{0,4})\geq k_D(w_{0,2},w_{0,3})-\frac{5}{4}-15C_0+\rho_3\\
	&\stackrel{\eqref{8-27-1}}{\geq} -\frac{5}{4}-15C_0-\rho_1+2\rho_3.
\end{aligned}
\ee
For any $v\in\gamma_{12}^{1}$, we know from Lemma \ref{9-27-1} that
\be\label{11-7}
\begin{aligned}k_D(v,w_{0,4})&\geq k_D(w_{0,4},w_{1,2})-k_D(w_{1,2},v)\geq  k_D(w_{0,4},w_{1,2})-\ell_k(\gamma_{12}^{0}[w_{1,2},v])\\ &\geq k_D(w_{0,4},w_{1,2})-\ell_k(\gamma_{12}^{0}[w_{1,2},w_{1,1}])
	\geq k_D(w_{0,4},w_{1,2})-k_D(w_{1,2},w_{1,1})-\frac{1}{2}\\
	&\stackrel{\eqref{8-26-1}+\eqref{23-08-20-2}}{\geq}-2-15C_0+\rho_2+\rho_3>\rho_3
\end{aligned}
\ee

Now, we see from \eqref{8-26-1} and \eqref{11-7} that $w_{1,1}$, $w_{1,2}$, $\gamma_{12}^1$  and $\varpi_1=w_{0,4}$  satisfy all the assumptions in Lemma \ref{lem-3.8}. It follows that
there are $w_{1,3}\in \gamma_{12}^{1}$ and $w_{1,4}\in \beta_0[w_{1,2},w_{0,4}]$ such that
the six-tuple $\Omega_1=[w_{1,1},w_{1,2},\gamma_{12}^1,w_{0,4},w_{1,3},w_{1,4}]$ satisfies Property \ref{Property-A} with constant $\rho_3$ $($see Figure \ref{fig-8-15-1}$)$.
Then Property \ref{Property-A}(\ref{3-3.2}) guarantees that for any $u\in \gamma_{12}^{1}$ and any $v\in \beta_0[w_{1,4},w_{0,4}]$,
$$k_D(u,v)\geq \rho_3-\frac{1}{4}.$$ This proves the claim.
\medskip

The construction of next six-tuple is formulated in the following claim.\medskip

\noindent {\bf Claim \ref{lem-24-05}-{\rm II}}: {\it Let $w_{2,1}=w_{1,3}$, $w_{2,2}=w_{1,4}$, $\gamma_{12}^2\in \Lambda_{w_{2,1}w_{2,2}}^\lambda$ and $\varpi_2=w_{0,4}$. Then   }
\begin{enumerate}
	\item\label{24-2-26-1}
	$k_D(w_{2,1},w_{2,2})=\rho_3$;
	\item\label{24-2-26-2}  {\it
		there are $w_{2,3}\in\gamma_{12}^{2}$
		and $w_{2,4}\in\beta_0[w_{2,2},\varpi_2]$ such that   }
	\begin{enumerate}
		\item\label{23-08-29-6}   {\it
			the six-tuple $\Omega_2^{\rho_3}=[w_{2,1},w_{2,2},\gamma_{12}^2,\varpi_2,w_{2,3},w_{2,4}]$ satisfies Property \ref{Property-A},   }
		\item\label{23-08-29-7}   {\it
			$k_D(w_{0,2},w_{2,2})=k_D(w_{1,2},w_{2,2})\geq \rho_3-\frac{1}{4}$, and    }
		\item\label{24-2-26-5}
		$k_D(w_{1,3},w_{2,3})=k_D(w_{2,1},w_{2,3})\leq \frac{1}{2}(\rho_1+\rho_2) +68A_0$.
\end{enumerate}\end{enumerate}
\smallskip 

\noindent {\bf Proof of Claim \ref{lem-24-05}-{\rm II}}. Since $\Omega_1^{\rho_3}$ satisfies Property \ref{Property-A}, we infer from Property \ref{Property-A}\eqref{3-3.1} that
\beqq\label{9-1-2}
k_D(w_{2,1},w_{2,2})=k_D(w_{1,3},w_{1,4})=\rho_3.
\eeqq This implies statement \eqref{24-2-26-1} of Claim \ref{lem-24-05}-{\rm II}.

Moreover, 
Lemma \ref{lem-22-09-30}$(2b)$ implies that for any $w\in\gamma_{12}^{2}$, it holds
\be\label{24-2-13-3}k_D(w,w_{0,4})>\rho_3.\ee
Thus, we know from statement  \eqref{24-2-26-1} of Claim \ref{lem-24-05}-{\rm II}, \eqref{24-2-13-3} and Lemma \ref{lem-3.8} (replacing $w_1,w_2,\varpi,\mathfrak{c}$ with $w_{2,1},w_{2,2},w_{1,4},\rho_3$) that there exist $w_{2,3}\in\gamma_{12}^{2}$ and $w_{2,4}\in \beta_0[w_{2,2},w_{0,4}]$, such that
the six-tuple $\Omega_2=[w_{2,1},w_{2,2},\gamma_{12}^2,w_{0,4},w_{2,3},w_{2,4}]$ satisfies Property \ref{Property-A} with constant $\rho_3$ $($see Figure \ref{fig-8-15-1}$)$. This proves the statement \eqref{23-08-29-6}.

Since $w_{0,2}=w_{1,2}\in \gamma_{12}^1$ and $w_{2,2}=w_{1,4}\in \beta_0[w_{1,2},\varpi_1]$, statement \eqref{23-08-29-7}
follows from Claim \ref{lem-24-05}-{\rm I}(\ref{23-08-29-5}). Furthermore, the statement \eqref{24-2-26-5} follows from Lemma \ref{lem-22-09-30} $(2)$, and hence, the proof of claim is complete.\medskip


The construction for the third six-tuple is given in the following claim.\medskip

\noindent {\bf Claim \ref{lem-24-05}-{\rm III}}: {\it Let $w_{3,1}=w_{2,3}$, $w_{3,2}=w_{2,4}$, $\gamma_{12}^3\in \Lambda_{w_{3,1}w_{3,2}}^\lambda$ and $\varpi_3=w_{0,4}$. Then   }
\begin{enumerate}
	\item\label{24-2-26-3}
	$k_D(w_{3,1},w_{3,2})=\rho_3$;
	\item\label{24-2-26-4}  {\it
		there are $w_{3,3}\in\gamma_{12}^{3}$
		and $w_{3,4}\in\beta_0[w_{3,2},\varpi_3]$ such that   }
	\begin{enumerate}
		\item\label{23-12-11-1}  {\it
			the six-tuple $\Omega_3^{\rho_3}=[w_{3,1},w_{3,2},\gamma_{12}^3,\varpi_3,w_{3,3},w_{3,4}]$ satisfies Property \ref{Property-A};    }
		\item\label{24-2-16-3}   {\it
			$k_D(w_{0,2},w_{3,2})=k_D(w_{1,2},w_{3,2})\geq \rho_3-\frac{1}{4}$, and  }
		\item\label{24-2-26-8}
		$k_D(w_{2,3},w_{3,3})=k_D(w_{3,1},w_{3,3})\leq \frac{1}{4}(\rho_1+\rho_2)+68A_0$.
\end{enumerate}\end{enumerate}
\smallskip 

\noindent {\bf Proof of Claim \ref{lem-24-05}-{\rm III}}. We have shown that both the six-tuples $\Omega_1^{\rho_3}$ and $\Omega_2^{\rho_3}$ satisfy Property \ref{Property-A}. Thus, statement \eqref{24-2-26-3} of Claim \ref{lem-24-05}-{\rm III} follows from Property \ref{Property-A}(\ref{3-3.1}).
Moreover, as Claim \ref{lem-24-05}-{\rm II}\eqref{24-2-26-5} shows
$$k_D(w_{2,1},w_{2,3})=k_D(w_{1,3},w_{2,3})\leq \frac{\rho_1+\rho_2}{2}+68A_0.$$
we infer from Lemma \ref{lem-22-09-30}$(3b)$ that for any $w\in\gamma_{12}^{3}$,
\be\label{24-5-20-1}
k_D(w,w_{1,4})>\rho_3.
\ee
Then Lemma \ref{lem-3.8} (replacing $w_1,w_2,\varpi,\mathfrak{c}$ with $w_{3,1},w_{3,2},w_{0,4},\rho_3$), together with the statement \eqref{24-2-26-3} of Claim \ref{lem-24-05}-{\rm III} and \eqref{24-5-20-1}, implies that there exist $w_{3,3}\in\gamma_{12}^{3}$ and $w_{3,4}\in \beta_0[w_{3,2},\varpi_3]$, such that
the six-tuple $\Omega_3=[w_{3,1},w_{3,2},\gamma_{12}^3,\varpi_3,w_{3,3},w_{3,4}]$ satisfies Property \ref{Property-A} with constant $\rho_3$ $($see Figure \ref{fig-8-15-1}$)$. This proves the statement \eqref{23-12-11-1}. Meanwhile, it follows from Lemma \ref{lem-22-09-30}$(2)$ and $w_{2,3}=w_{3,1}$ that statement \eqref{24-2-26-8} holds. Hence, the proof of claim is complete.


\medskip

Finally, we construct the last six-tuple in the following claim.\medskip

\noindent {\bf Claim \ref{lem-24-05}-{\rm IV}}: {\it Let $w_{4,1}=w_{3,3}$, $w_{4,2}=w_{3,4}$, $\gamma_{12}^4\in \Lambda_{w_{4,1}w_{4,2}}^\lambda$ and $\varpi_4=w_{0,4}$. Then   }
\begin{enumerate}
	\item\label{24-2-26-3 IV}
	$k_D(w_{4,1},w_{4,2})=\rho_3$;
	\item\label{24-2-26-4 IV}  {\it
		there are $w_{4,3}\in\gamma_{12}^{4}$
		and $w_{4,4}\in\beta_0[w_{4,2},\varpi_4]$ such that   }
	\begin{enumerate}
		\item\label{23-12-11-1 IV}  {\it
			the six-tuple $\Omega_4=[w_{4,1},w_{4,2},\gamma_{12}^4,\varpi_4,w_{4,3},w_{4,4}]$ satisfies Property \ref{Property-A} with constant $\rho_3$;    }
		\item\label{24-2-16-3 IV}   {\it
			$k_D(w_{0,2},w_{4,2})=k_D(w_{1,2},w_{4,2})\geq \rho_3-\frac{1}{4}$, and  }
		\item\label{24-2-26-8 IV}
		$k_D(w_{3,3},w_{4,3})=k_D(w_{4,1},w_{4,3})\leq  \frac{1}{8}(\rho_1+\rho_2)+68A_0$.
	\end{enumerate}
\end{enumerate}
\smallskip

\noindent {\bf Proof of Claim \ref{lem-24-05}-{\rm IV}}. Since all the six-tuples $\Omega_1^{\rho_3},$ $\Omega_2^{\rho_3}$ and $\Omega_3^{\rho_3}$ satisfy Property \ref{Property-A}, statement \eqref{24-2-26-3 IV} of Claim \ref{lem-24-05}-{\rm IV} follows from Property \ref{Property-A}(\ref{3-3.1}).
Moreover, as Claim \ref{lem-24-05}-{\rm III}\eqref{24-2-26-5} shows
$$k_D(w_{3,1},w_{3,3})=k_D(w_{2,3},w_{3,3})\leq \frac{1}{4}(\rho_1+\rho_2) +68A_0,$$
which, together with {\rm Claim \ref{lem-24-05}-}{\rm II}, gives that
$$k_D(w_{1,3},w_{3,3})\leq k_D(w_{1,3},w_{2,3})+k_D(w_{2,3},w_{3,3})\leq \frac{3}{4}(\rho_1+\rho_2)+136A_0.$$
It follows from Lemma \ref{lem-22-09-30}$(3b)$ that for any $w\in\gamma_{12}^{4}$,
\be\label{H25-0812-1}
k_D(w,w_{0,4})>\rho_3.
\ee
Applying Lemma \ref{lem-3.8} (replacing $w_1,w_2,\varpi,\mathfrak{c}$ with $w_{4,1},w_{4,2},w_{0,4},\tau_3$), together with the statement \eqref{24-2-26-3 IV} of Claim \ref{lem-24-05}-{\rm IV} and \eqref{H25-0812-1}, gives that there exist $w_{4,3}\in\gamma_{12}^{4}$ and $w_{4,4}\in \beta_0[w_{4,2},\varpi_4]$, such that
the six-tuple $\Omega_4=[w_{4,1},w_{4,2},\gamma_{12}^4,\varpi_4,w_{4,3},w_{4,4}]$ satisfies Property \ref{Property-A} with constant $\rho_3$ $($see Figure \ref{fig-8-15-1}$)$. This proves the statement \eqref{23-12-11-1 IV}. Since $w_{3,3}=w_{4,1}$, the statement \eqref{24-2-26-8 IV} holds by Lemma \ref{lem-22-09-30} $(2)$, and hence, the proof of claim is complete.
\medskip

\epf


\subsection{Proofs of the main results}\label{subsec:final section}
In this section, we prove our main results.	
Set
\be\label{24-2-28-2}
\tau_2=\frac{1}{64}\log \kappa\;\;\mbox{and}\;\;\tau_3=\frac{1}{16}\log \kappa.
\ee 

For $i\in \{1, \ldots, m\}$, let $\xi_i$ be the points defined in Section \ref{sec-3}, and $w_{0,2}=\xi_1$ (see Figure \ref{fig3}).
Since Lemma \ref{cl3.2-1}$(i)$ and the choice of $\xi_i$ in \eqref{2022-09-30-3} ensure that
\be\label{23-08-20-1}
k_D(\xi_{i},y_0)= k_D(w_{0,2},y_0)\geq \log \frac{\kappa-2}{2{\color{blue}c_2}(2+\kappa_1)}\stackrel{\eqref{24-2-28-2}}{\geq}2\tau_3+\frac{1}{2},
\ee
we infer from the above estimate that there exists $w_{0,1}\in\gamma_1[w_{0,2},y_0]$ such that
\be\label{2022-10-12-2}
k_D(w_{0,1},w_{0,2})=\tau_3.
\ee 

Let  $\gamma_{12}^0=\gamma_1[w_{0,1},w_{0,2}]$. Then for any $w\in\gamma_{12}^0$, we know from Lemma \ref{9-27-1} that
\beqq
k_D(w,\xi_{m}) \geq 
\ell_k(\gamma_1[w,\xi_{m}])-\frac{1}{2} >\ell_k(\gamma_1[\xi_{m},y_0])-\frac{1}{2}
\geq k_D(\xi_{m},y_0)-\frac{1}{2}\stackrel{\eqref{23-08-20-1}}{>}2\tau_3.
\eeqq
Thus it follows from \eqref{2022-10-12-2} and Lemma \ref{lem-3.8} that there exist $w_{0,3}\in \gamma_{12}^0$ and $w_{0,4}\in \beta_0[w_{0,2},\xi_{m}]$ such that for $\varpi=\xi_{m}$, the six-tuple $[w_{0,1},w_{0,2},\gamma_{12}^0,\varpi,w_{0,3},w_{0,4}]$ satisfies Property \ref{Property-A} with constant $\kappa_3$.
Then by Lemma \ref{lem-23-01}, we have
\be\label{24-09-14-1}
k_D(w_{0,1},w_{0,3})\leq 11A_0.
\ee

Next, we let $w_{1,1}^0=w_{0,1}$, $w_{1,2}^0=w_{0,2}$, $w_{1,3}^0=w_{0,3}$, $w_{1,4}^0=w_{0,4}$,
$\varpi^0=\varpi$, $\gamma^0=\gamma_{12}^0$ and $\Omega^0=[w_{1,1}^0, w_{1,2}^0, \gamma^0, \varpi^0, w_{1,3}^0, w_{1,4}^0]$. Take $C_1=[\ell_k(\beta_0)]$, $\tau_{1,0}=\tau_1=11A_0$ and $\tau_{2,0}=\tau_2$. For each $i\in\{1,\cdots, C_1\}$, set $\tau_{1,i}=\frac{\tau_{1,i-1}+\tau_{2,i-1}}{2^3}+68A_0$ and $\tau_{2,i}=\tau_2+\tau_{1,i}$.

Then we shall prove the following lemma, which is the key for the proof of Theorem \ref{thm:main}.

\begin{lem}\label{lem-24}
	There exist $(C_1+1)$ six-tuples $\Omega_i^{\tau_3}=[w_{1,1}^i, w_{1,2}^i, \gamma^i, \varpi^i, w_{1,3}^i, w_{1,4}^i]$ with Property \ref{Property-A} such that the following conclusions hold:
	\ben
	\item\label{H25-0813-1} For each $i\in\{1,\cdots,C_1\}$, $w_{1,2}^i\in\beta_0[w_{1,2}^{i-1},w_{1,4}^{i-1}]$, $\varpi^i=w_{1,4}^{i-1}$ and $\gamma^i=\gamma_{w_{1,1}^iw_{1,2}^i}^{\lambda}$ with $\gamma_{w_{1,1}^iw_{1,2}^i}^{\lambda} \in \Lambda_{w_{1,1}^iw_{1,2}^i}$.
	
	\item\label{H25-0813-2} For each $i\in\{1,\cdots,C_1\}$, $k_D(w_{1,1}^i,w_{1,3}^i)\leq \tau_{1,i}$.
	
	\item\label{H25-0813-3} For each $i\in\{1,\cdots,C_1\}$, $k_D(w_{1,2}^{i-1},w_{1,2}^i)>\tau_3-\frac{1}{4}$.
	\een
\end{lem}
\bpf Since $\tau_{2,0}>5\tau_{1,0}+24A_0e^{36A_0}$, for $i=1$, the result follows from Lemma \ref{lem-24-05}, \eqref{2022-10-12-2} and \eqref{24-09-14-1}. The general case can be proved by a standard induction argument, for which we formulate as a claim below.
\bcl\label{H25-0814}
For $k\geq 1$, if Lemma \ref{lem-24-05} holds for each $i\in\{1,\cdots,k\}$, then it holds for $i=k+1$ as well.
\ecl

An elementary calculation yields
$$\tau_{1,k}=68A_0+\frac{11A_0}{2^{1+2k}}+\sum_{i=1}^k \frac{68A_0}{2^{2i-1}}+\sum_{i=1}^k \frac{1}{2^{2i+1}}\tau_2,$$
which, together with the fact ``$\tau_{2,k}=\tau_2+\tau_{1,k}$", implies
$$\tau_3>2\tau_2>\tau_{2,k}>5\tau_{1,k}+24A_0e^{6A_0}.$$
Thus, we may apply Lemma \ref{lem-24-05} again to conclude that Lemma \ref{lem-24} holds for $i=k+1$. 

Lemma \ref{lem-24} follows immediately from Claim \ref{H25-0814}.
\epf

\begin{proof}[Proof of Theorem \ref{thm:main}]
	By Lemma \ref{lem-24}(\ref{H25-0813-3}), we have
	$$\ell_k(\beta_0)\geq \sum_{i=2}^{C_1}k_D(w_{1,2}^{i-1},w_{1,2}^i)\geq C_1(\tau_3-\frac{1}{4})>\ell_k(\beta_0).$$
	This contradiction shows that \eqref{lem2-eq-2} cannot hold. The proof of Theorem \ref{thm:main} is thus complete.
\end{proof}

\begin{proof}[Proof of Corollary \ref{coro1}]
	Since the Gromov hyperbolicity is preserved under $(M,C)$-CQH by \cite[Theorem 3.18]{Vai10}, Corollary \ref{coro1} follows directly from Theorem \ref{thm:main}.  
\end{proof}


\end{document}